\documentclass[twoside,11pt]{article}

\usepackage{amsmath,amsfonts}

\newcommand{\D}{\displaystyle}

\def \qed {\hfill \vrule height6pt width 6pt depth 0pt}

\begin{document}
\title{{\bf {Decay of the solution to the bipolar Euler-Poisson system with damping in $\mathbb{R}^3$}}}
\author{ Z{\sc higang}
W{\sc u}\thanks{Corresponding author. E-mail: mathzgwu@yahoo.com.cn}\\
\small{\it Department of Mathematics, Hangzhou Normal
University, Hangzhou, P.R. China}\\
 \and  W{\sc eike} W{\sc ang}\\
{\small{\it Department of Mathematics, Shanghai Jiao Tong
University, Shanghai, P.R. China}} }

\date{}
\maketitle

\textbf{{\bf Abstract:}} We construct the global solution to the
Cauchy's problem of the bipolar Euler-Poisson equations with damping
in $\mathbb{R}^3$ when $H^3$ norm of the initial data is small. If
further, the $\dot{H}^{-s}$ norm ($0\leq s<3/2)$ or
$\dot{B}_{2,\infty}^{-s}$ norm ($0<s\leq3/2$) of the initial data is
bounded, we give the optimal decay rates of the solution. As a
byproduct,  the decay results of the $L^p-L^2$ ($1\leq p\leq2$) type
hold without the smallness of the $L^p$ norm of the initial data. In
particular, we deduce that $\|\nabla^k(\rho_1-\rho_2)\|_{L^2}
\sim(1+t)^{-\frac{5}{4}-\frac{k}{2}}$ and
$\|\nabla^k(\rho_i-\bar{\rho},u_i,\nabla\phi)\|_{L^2}
\sim(1+t)^{-\frac{3}{4}-\frac{k}{2}}$. We improve the decay results
in Li and Yang \cite{Li3}(\emph{J.Differential Equations} 252(2012),
768-791), where they showed the decay rates as
$\|\nabla^k(\rho_i-\bar{\rho})\|_{L^2}
\sim(1+t)^{-\frac{3}{4}-\frac{k}{2}}$ and
$\|\nabla^k(u_i,\nabla\phi)\|_{L^2}
\sim(1+t)^{-\frac{1}{4}-\frac{k}{2}}$, when the $H^3\cap L^1$ norm
of the initial data is small. Our analysis is motivated by the
technique developed recently in Guo and Wang \cite{Guo}(\emph{Comm.
Partial Differential Equations} 37(2012), 2165-2208) with some
modifications.

 \bigbreak \textbf{{\bf Key Words}:}
Bipolar Euler-Poisson system; Global existence; Decay estimates;
Negative Sobolev's space; Negative Besov's space.

\bigbreak
\section*{ 1.\ Introduction}

\quad\quad We consider the compressible bipolar Euler-Poisson
equations with damping (BEP)
$$
\left\{\begin{array}{l}
\partial_t\rho_1+{\rm div} (\rho_1 u_1)=0, \\[2mm]
\partial_t(\rho_1 u_1)+ {\rm div}(\rho_1u_1\otimes u_1)+\nabla
P(\rho_1)
=\rho_1\nabla\phi-\rho_1u_1, \\[2mm]
\partial_t\rho_2+{\rm div} (\rho_2 u_2)=0, \\[2mm]
\partial_t(\rho_2 u_2)+ {\rm div}(\rho_2u_2\otimes u_2)+\nabla
P(\rho_2)
=-\rho_2\nabla\phi-\rho_2u_2, \\[2mm]
\Delta\phi=\rho_1-\rho_2,\ \ \ x\in\mathbb{R}^3,\ t\geq0,
\end{array}
        \right.
        \eqno({\rm 1.1})
$$
where the unknown functions $\rho_i(x,t),u_i(x,t)\
(i=1,2),\phi(x,t)$ represent the charge densities, current
densities, velocities and electrostatic potential, respectively, and
the pressures $P=P(\rho_i)$ is a smooth function with $P'(\rho_i)>0$
for $\rho_i>0$. The system (1.1) is usually described charged
particle fluids, for example, electrons and holes in semiconductor
devices, positively and negatively charged ions in a plasma. We
refer to \cite{Jugel,Sitenko} for the physical background of the
system (1.1).

In this paper, we will study the global existence and large time
behavior of the smooth solutions for the system (1.1) with the
following initial data
$$
\rho_i(x,0)=\rho_{i0}(x)>0,\ u_i(x,0)=u_{i0}(x),\ i=1,2. \eqno(1.2)
$$

A lot of important works have been done on system (1.1). For
one-dimensional case, we refer to Zhou and Li \cite{Zhou} and Tsuge
\cite{Tsuge} for the unique existence of the stationary solutions,
Natalini \cite{Natalini1} and Hsiao and Zhang \cite{Hsiao1} for
global entropy weak solutions in the framework of compensated
compactness on the whole real line and bounded domain respectively,
Natalini \cite{Natalini1}, Hsiao and Zhang \cite{Hsiao2} for the
relaxation-time limit, Gasser and Marcati \cite{Gasser2} for the
combined limit, Huang and Li \cite{Huang} for the large-time
behavior and quasi-neutral limit of $L^\infty$-solution, Zhu and
Hattori \cite{Zhu} for the stability of steady-state solutions to a
recombined one-dimensional bipolar hydrodynamical model, Gasser,
Hsiao and Li \cite{Gasser1} for large-time behavior of smooth small
solution.

For the multi-dimensional case, Lattanzio \cite{Lattanzio} discussed
the relaxation limit, and Li \cite{Li2} considered the diffusive
relaxation. Ali and J\"{u}ngel \cite{Ali} and Li and Zhang
\cite{Li1} studied the global smooth solutions of the Cauchy problem
in the Sobolev's space and Besov's space, respectively. Later, Ju
\cite{Ju} investigated the global existence of smooth solution to
the IBVP for the 3D bipolar Euler-Poisson system (1.1).

Recently, Using the classical energy method together with the
analysis of the Green's function, Li and Yang \cite{Li3}
investigated the optimal decay rate of the Cauchy's problem of the
system (1.1) of the classical solution when the initial data is
small in the space $H^3\cap L^1$. They deduced that the electric
field (a nonlocal term in Hyperbolic-parabolic system) slows down
the decay rate of the velocity of the BEP system, also see the
recent works \cite{LiH1,LiH2,Zhang,Wang1,Wu1,Wu2} on the decay of
the solutions to the unipolar Navier-Stokes-Poisson equations (NSP)
and unipolar Euler-Poisson equations with damping. In fact, by the
detailed analysis of the Green's function, it shows that the
presence of the electric field field slows down the decay rate in
$L^2$-norm of the velocity of the unipolar NSP system with the
factor 1/2 comparing with the Navier-Stokes system (NS) when the
initial perturbation $\rho_0-\bar{\rho},u_0\in L^p\cap H^3$ with
$p\in[1,2]$.

However, Wang \cite{WangY1} gave a different comprehension of the
effect of the electric field on the time decay rates of the solution
of the unipolar NSP system. The key idea is making an instead
assumption on the initial perturbation $\rho_0-\bar{\rho}\in
\dot{H}^{-1},u_0\in L^2$. As a result, the electric field does not
slow down but rather enhances the time decay rate of the density
with the factor $\frac{1}{2}$. The method in \cite{WangY1} is
initially established in Guo and Wang \cite{Guo} for the estimates
in the negative Sobolev's space. The proof in \cite{Guo} is based on
a family of energy estimates with minimum derivative counts and
interpolations among them without linear decay analysis. Very
recently, using this kind of energy estimates, Tan and Wang
\cite{Tan} discussed the Euler equations with damping in
$\mathbb{R}^3$, where they also gave the estimates in the negative
Besov's space.

The main purpose of this paper is to improve the $L^2$-norm decay
estimates of the solutions in Li and Wang \cite{Li3} by using this
refined energy method together with the interpolation trick in
\cite{Guo,WangY1,Tan}. Comparing with \cite{Guo,WangY1,Tan}, the
main additional difficulties are due to the presence of electronic
field and the couple of two carriers by the Poisson equation. First,
as Wang \cite{WangY1} pointed out, for the bipolar NSP system, there
is one term $n_iu_i\nabla\phi$ can not be controlled by the
dissipation terms when using this refined energy method, see the
details in \cite{WangY1}. However, after an elaborate calculation,
we can get each $l$-th ($l=0,1,2,3$) level energy estimate for the
BEP system (1.1), see (2.26)-(2.28) and (2.38)-(2.39) in Lemma 2.10
and Lemma 2.11. Second, one can not obtain the dissipation term for
$\|\rho_i\|_{L^2}$ in the energy estimates as the unipolar case in
\cite{WangY1} since two species are strongly coupled by the Poisson
equation for bipolar case. In fact, we only can get the estimate
$\|\nabla^k(\rho_1-\rho_2)\|_{L^2}\leq\|\nabla^{k+1}\nabla\phi\|_{L^2}$
for the BEP system (1.1). As a result, one can not directly deal
with the case $s\in(\frac{1}{2},\frac{3}{2})$ for the estimates in
the negative Sobolev's space or negative Besov's space by using the
decay result for the case $s\in[0,\frac{1}{2}]$ as in \cite{WangY1}.
In fact, the proof for the case $s\in(\frac{1}{2},\frac{3}{2})$ in
Wang \cite{WangY1} strongly depends on the derived decay result of
the case $s=\frac{1}{2}$:
$\|\rho\|_{L^2}\leq\|\nabla\nabla\phi\|_{L^2}\sim(1+t)^{-\frac{l+s}{2}}=(1+t)^{-\frac{3}{4}}$.
After a detailed analysis, and by separating the cases that
$s\in[0,\frac{1}{2}]$, $s\in(\frac{1}{2},1)$ and
$s\in[1,\frac{3}{2})$ for the space $\dot{H}^{-s}$ and
$s\in[0,\frac{1}{2}]$, $s\in(\frac{1}{2},1)$, $s\in[1,\frac{3}{2})$
and $s=\frac{3}{2}$ for the space $\dot{B}_{2,\infty}^{-s}$, we
achieve these estimates (See Lemma 2.13, Lemma 2.14 and Subsection
3.2).

Our main results are stated in the following theorems:

\bigbreak \noindent\textbf{Theorem 1.1.} Let $P'(\rho_i)>0 (i=1,2)$
for $\rho_i>0$, and $\bar{\rho}>0$. Assume that
$(\rho_i-\bar{\rho},u_{i0},\nabla\phi_0)\in H^3(\mathbb{R}^3)$ for
$i=1,2$, with
$\epsilon_0=:\|(\rho_{i0}-\bar{\rho},u_{i0},\nabla\phi_0)\|_{H^3(\mathbb{R}^3)}$
small. Then there exists a unique, global, classical solution
$(\rho_1,u_1,\rho_2,u_2,\phi)$ satisfying that for all $t\geq0$,
$$\arraycolsep=1.5pt  \begin{array}[b]{rl}
&\D\|(n_1,u_1,n_2,u_2,\nabla\phi)\|_{H^3}^2+\int_0^t\|(u_1,u_2)\|_{H^3}^2+\|(\nabla
n_1,\nabla n_2,\nabla(\nabla\phi))\|_{H^2}^2d\tau\\[3mm]
\leq &\D C\|(u_{10},,u_{20},n_{10},n_{20},\nabla\phi_0)\|_{H^3}^2.
\end{array}
\eqno(1.3)
$$

\bigbreak \noindent\textbf{Remark 1.1.} From the fact
$\nabla\phi_0\in H^3$ is equivalent to $(n_{10}-n_{20})\in
H^2\cap\dot{H}^{-1}$ deriving from the poisson equation $(1.1)_5$,
we can replace the initial assumption $\nabla\phi_0\in H^3$ by
$n_{i0}\in H^2\cap H^{-1}$.

\bigbreak \noindent\textbf{Theorem 1.2.} Under the assumptions of
Theorem 1.1. If further,
$(\rho_{i0}-\bar{\rho},u_{i0},\nabla\phi_0)\in \dot{H}^{-s}$ for
some $s\in[0,3/2)$ or
$(\rho_{i0}-\bar{\rho},u_{i0},\nabla\phi_0)\in\dot{B}_{2,\infty}^{-s}$
for some $s\in(0,3/2]$, then for all $t\geq0$, there exists a
positive constant $C_0$ such that
$$
\|(\rho_i-\bar{\rho},u_i,\nabla\phi)(t)\|_{\dot{H}^{-s}}\leq
C_0\eqno(1.4)
$$
or
$$
\|(\rho_i-\bar{\rho},u_i,\nabla\phi)(t)\|_{\dot{B}_{2,\infty}^{-s}}\leq
C_0,\eqno(1.5)
$$
and
$$
\|\nabla^l(\rho_i-\bar{\rho},u_i,\nabla\phi)(t)\|_{H^{3-l}}\leq
C_0(1+t)^{-\frac{l+s}{2}}\ {\rm for}\ l=0,1,2;\
s\in[0,\frac{3}{2}];\eqno(1.6)
$$
$$
\|\nabla^l(\rho_1-\rho_2)(t)\|_{L^2}\leq
C_0(1+t)^{-\frac{l+s+1}{2}}\ {\rm for}\ l=0,1;\ s\in[0,\frac{3}{2}].
\eqno(1.7)
$$

\bigbreak \noindent\textbf{Remark 1.2.} (1.7) is derived from (1.6)
and the fact
$$
\|\nabla^l(\rho_1-\rho_2)\|_{L^2}=\|\nabla^l\Delta\phi\|_{L^2}\leq\|\nabla^{l+1}\nabla\phi\|_{L^2},
$$
which shows the presence of the electric field enhances the time
decay rate of disparity between two species.

Note that Lemma 2.4 (the Hardy-Littlewood-Sobolev theorem) implies
that for $p\in(1,2],L^p\subset\dot{H}^{-s}$ with
$s=3(\frac{1}{p}-\frac{1}{2})$ and Lemma 2.6 implies that for
$p\in[1,2),L^p\subset\dot{B}_{2,\infty}^{-s}$ with
$s=3(\frac{1}{p}-\frac{1}{2})$. Then Theorem 1.2 yields the
following usual $L^p-L^2$ type of optimal decay results.

\bigbreak \noindent\textbf{Corollary 1.1.} Under the assumptions of
Theorem 1.2 except that we replace the $\dot{H}^{-s}$ or
$\dot{B}_{2,\infty}^{-s}$ assumption by that
$(\rho_{i0}-\bar{\rho},u_{i0},\nabla\phi_0)\in L^p$ for some
$p\in[1,2]$, then the following decay results hold:
$$
\|\nabla^l(\rho_i-\bar{\rho},u_i,\nabla\phi)(t)\|_{H^{3-l}}\leq
C_0(1+t)^{-\frac{3}{2}(\frac{1}{p}-\frac{1}{2})-\frac{l}{2}},\ {\rm
for}\ l=0,1,2; \eqno(1.8)
$$
$$
\|\nabla^l(\rho_1-\rho_2)(t)\|_{L^2}\leq
C_0(1+t)^{-\frac{3}{2}(\frac{1}{p}-\frac{1}{2})-\frac{l+1}{2}},\
{\rm for}\ l=0,1.\eqno(1.9)
$$

\bigbreak \noindent\textbf{Remark 1.3.} From Corollary 1.1, we know
the each order derivatives of the density $\rho_i-\bar{\rho}$ and
velocity $u_i$ have the same decay rate in $H^2$ norm as the
solutions of the Navier-Stokes equations. While the velocity $u_i$
in \cite{Li3} decays at the rate $(1+t)^{-\frac{1}{4}}$ in $L^2$
norm which is slower than the rate $(1+t)^{-\frac{3}{4}}$ for the
compressible Navier-Stokes equations.

\bigbreak \noindent\textbf{Remark 1.4.} The energy method (close the
energy estimates at each $l$-th level with respect to the spatial
derivatives of the solutions) in this paper can not be applied to
the bipolar Navier-Stokes-Poisson equations. In fact, as Wang
\cite{WangY1} pointed out, there is one term $n_iu_i\nabla\phi$ can
not be controlled by the dissipation terms, see the Introduction in
\cite{WangY1}. Hence, it is also interesting to apply this energy
method to the bipolar Navier-Stokes-Poisson equations.

\bigbreak \noindent\textbf{Remark 1.5.} We also notice that, the
similar arguments can be used to investigate the full
(nonisentropic) bipolar hydrodynamic models, which is under
consideration.

\bigbreak \noindent\textbf{Notations.} In this paper, $\nabla^l$
with an integer $l\geq 0$ stands for the usual any spatial
derivatives of order $l$. For $1\leq p\leq \infty$ and an integer
$m\geq 0$, we use $L^p$ and $W^{m,p}$ denote the usual Lebesgue
space $L^p(\mathbb{R}^n)$ and Sobolev spaces $W^{m,p}(\mathbb{R}^n)$
with norms $\|\cdot\|_{L^p}$ and $\|\cdot\|_{W^{m,p}}$,
respectively, and set $H^m=W^{m,2}$ with norm $\|\cdot\|_{H^m}$ when
$p=2$. In addition, for $s\in\mathbb{R}$, we define a
pseudo-differential operator $\Lambda^s$ by
$$
\Lambda^sg(x)=\int_{\mathbb{R}^n}|\xi|^s\hat{g}(\xi){\rm
e}^{2\pi\sqrt{-1}x\cdot\xi} {\rm d}\xi,
$$
where $\hat{g}$ denotes the Fourier transform of $g$. We define the
homogeneous Sobolev's space $\dot{H}^s$ of all $g$ for which
$\|g\|_{\dot{H}^s}$ is finite, where
$$
\|g\|_{\dot{H}^s}:=\|\Lambda^sg\|_{L^2}=\||\xi|^s\hat{g}\|_{L^2}.
$$
Let $\eta\in C_0^\infty(\mathbb{R}_\xi^3)$ be such that
$\eta(\xi)=1$ when $|\xi|\leq1$ and $\eta(\xi)=0$ when $\xi\geq2$.
We define the homogeneous Besov's spaces
$\dot{B}_{2,\infty}^{-s}(\mathbb{R}^3$ with norm
$\|\cdot\|_{\dot{B}_{p,r}^{-s}}$ defined by
$$
\|f\|_{\dot{B}_{p,r}^{-s}}:=(\sum_{j\in
\mathbb{Z}}2^{rsj}\|\dot{\Delta}_jf\|_{L^p}^r)^{\frac{1}{r}}.
$$
Here $\dot{\Delta}_{j}f:=F^{-1}(\varphi_j)\ast f$,
$\varphi(\xi)=\eta(\xi)-\eta(2\xi)$ and
$\varphi_j(\xi)=\varphi(2^{-j}\xi)$.

Throughout this paper, we will use a non-positive index $s$. For
convenience, we will change the index to be ``$-s$'' with $s\geq0$.
$C$ or $C_i$ denotes a positive generic (generally large) constant
that may vary at different places. For simplicity, we write $\int
f:=\int_{\mathbb{R}^3}f{\rm d}x$.

The rest of the paper is arranged as follows. In section 2, we give
some useful Sobolev's inequality and Besov's inequality, then we
give energy estimate in $H^3$ norm and some estimates in
$\dot{H}^{-s}$ and $\dot{B}_{2,\infty}^{-s}$. The proof of global
existence and temporal decay results of the solutions will be
derived in Section 3.

\section*{ 2.\ Nonlinear energy estimates}

\subsection*{ 2.1.\ Preliminaries}

\quad\quad In this subsection we give some Sobolev's inequalities
and Besov's inequalities, which will be used in the next sections.

\bigbreak \noindent\textbf{Lemma 2.1.} (Gagliardo-Nirenberg's
inequality). Let $0\leq m,k\leq l$, then we have
$$
\|\nabla^k g\|_{L^p}\leq C \|\nabla^m
g\|_{L^q}^{1-\theta}\|\nabla^lg\|_{L^r}^\theta,
$$
where $k$ satisfies
$$
\frac{1}{p}-\frac{k}{n}=(1-\theta)\left(\frac{1}{q}-\frac{m}{n}\right)+\theta\left(\frac{1}{r}-\frac{l}{n}\right).
$$

\bigbreak \noindent\textbf{Lemma 2.2.} (Moser-type calculus) (i) Let
$k\geq1$ be an integer and define the commutator
$$
[\nabla^k,g]h=\nabla^k(gh)-g\nabla^kh.
$$
Then we have
$$
\|[\nabla^k,g]h\|_{L^2}\leq C_k(\|\nabla
g\|_{L^\infty}\|\nabla^{k-1}h\|_{L^2}+\|\nabla^kg\|_{L^2}\|h\|_{L^\infty}).
$$
(ii) If $F(\cdot)$ is a smooth function, $f(x)\in H^k\cap L^\infty$,
then we have
$$
\|\nabla^k F(f)\|\leq C(k,F,\|f\|_{L^\infty})\|\nabla^k f\|.
$$

\bigbreak \noindent\textbf{Lemma 2.3.} (\cite{Guo}, Lemma A.5) Let
$s\geq0$ and $l\geq0$, then we have
$$
\|\nabla^lg\|_{L^2}\leq
C\|\nabla^{l+1}g\|_{L^2}^{1-\theta}\|g\|_{\dot{H}^{-s}}^\theta,\
{\rm where}\ \theta=\frac{1}{l+s+1}.
$$

\bigbreak \noindent\textbf{Lemma 2.4.} (\cite{Stein}, Chapter V,
Theorem 1) Let $0<s<n, 1<p<q<\infty,
\frac{1}{q}+\frac{s}{n}=\frac{1}{p}$, then
$$
\|\Lambda^{-s}g\|_{L^q}\leq C\|g\|_{L^p}.
$$

Next, we give some lemmas on Besov space $\dot{B}_{2,\infty}^{-s}$.

\bigbreak \noindent\textbf{Lemma 2.5.} (\cite{Tan}) Suppose $k\geq0$
and $s>0$, then we have
$$
\|\nabla^kf\|_{L^2}\leq
C\|\nabla^{k+1}f\|_{L^2}^{1-\theta}\|f\|_{\dot{B}_{2,\infty}^{-s}}^\theta,\
{\rm where}\ \theta=\frac{1}{l+1+s}.
$$

\bigbreak \noindent\textbf{Lemma 2.6.} (\cite{Sohinger}) Suppose
that $s>0$ and $1\leq p<2$. We have the embedding
$L^p\subset\dot{B}_{q,\infty}^{-s}$ with $1/2+s/3=1/p$. In
particular we have the estimate
$$
\|f\|_{\dot{B}_{2,\infty}^{-s}}\leq C\|f\|_{L^p}.
$$

\bigbreak \noindent\textbf{Lemma 2.7.} (\cite{Sohinger}) Suppose
$k\geq0$ and $s>0$, then we have
$$
\|\nabla^k\|_{L^2}\leq
C\|\nabla^{k+1}f\|_{L^2}^{1-\theta}\|f\|_{\dot{B}_{2,\infty}^{-s}}^\theta,\
{\rm where}\ \theta=\frac{1}{l+s+1}.
$$

\bigbreak \noindent\textbf{Lemma 2.8.} (\cite{Stein}) If $1\leq
r_1\leq r_2\leq\infty$, then
$$
\dot{B}_{2,r_1}^{-s}\in\dot{B}_{2,r_2}^{-s}.
$$

\bigbreak \noindent\textbf{Lemma 2.9.} (\cite{Sohinger}) If $m>l\geq
k$ and $1\leq p\leq q\leq r\leq\infty$. We have
$$
\|g\|_{\dot{B}_{2,q}^{l}}\leq
C\|g\|_{\dot{B}_{2,r}^{k}}^\theta\|g\|_{\dot{B}_{2p}^{m}}^{1-\theta},
$$
where $l=k\theta+m(1-\theta),\
\frac{1}{q}=\frac{\theta}{r}+\frac{1-\theta}{p}$.

\subsection*{ 2.2\ Energy estimates in $H^3$-norm}

\quad\quad We reformulate the nonlinear system (1.1) for
$(\rho_1,u_1,\rho_2,u_2)$ around the equilibrium state
$(\bar{\rho},0,\bar{\rho},0)$. Without loss of generality, we can
assume $\bar{\rho}=1$ and $P'(\bar{\rho})=1$. Denote
$$
n_i=\rho_i-1,\ h(n_i)=\frac{P'(\rho_i)}{\rho_i}-1,
$$
then the Cauchy problem for $(n_1,u_1,n_2,u_2,\phi)$ is given by
$$
\left\{\begin{array}{l}
\partial_t n_1+{\rm div} u_1=-u_1\cdot\nabla n_1-n_1{\rm div}u_1, \\[2mm]
\partial_t u_1+u_1+\nabla n_1-\nabla\phi=-u_1\cdot\nabla
u_1-h(n_1)\nabla n_1, \\[2mm]
\partial_t n_2+{\rm div} u_2=-u_2\cdot\nabla n_2-n_2{\rm div}u_2, \\[2mm]
\partial_t u_2+u_2+\nabla n_2+\nabla\phi=-u_2\cdot\nabla
u_2-h(n_2)\nabla n_2, \\[2mm]
\Delta\phi=n_1-n_2,\\[2mm]
(n_1,u_1,n_2,u_2)(x,0)=(\rho_{10}-1,u_{10},\rho_{20}-1,u_{20})(x).
\end{array}
        \right.
        \eqno({\rm 2.1})
$$

In this section, we will derive a priori nonlinear energy estimates
for the equivalent system (2.1). Hence we assume a priori assumption
that for a sufficiently small constant $\delta>0$,
$$
\|n_i(t)\|_{H^3}+\|u_i(t)\|_{H^3}+\|\nabla\phi(t)\|_{H^3}\leq\delta,\
i=1,2, \eqno(2.2)
$$
which together with Sobolev's inequality, we have the facts
$$
1/2\leq n_i\leq2,\ |h^{(k)}(n_i)|\leq C,\ i=1,2,\ {\rm for\ any\
k\geq0}.\eqno(2.3)
$$

We first deduce the following energy estimates which contains the
dissipation estimate for $u_1,u_2$.

\bigbreak \noindent\textbf{Lemma 2.10.} Assume that $0\leq k\leq 2$,
then we have
$$\arraycolsep=1.5pt  \begin{array}[b]{rl}
&\D\frac{1}{2}\frac{d}{dt}\int|\nabla^k(n_1,u_1,n_2,u_2,\nabla\phi)|^2+\|\nabla^k(u_1,u_2)\|_{L^2}^2\\[3mm]
\leq &\D
C\delta(\|\nabla^{k+1}n_1\|_{L^2}^2\!+\!\|\nabla^ku_1\|_{L^2}^2\!+\!\|\nabla^{k+1}n_2\|_{L^2}^2
\!+\!\|\nabla^ku_2\|_{L^2}^2\!+\!\|\nabla^{k+1}\nabla\phi\|_{L^2}^2).
\end{array}
\eqno(2.4)
$$
\noindent\textbf{Proof.} For $0\leq k\leq 2$, applying $\nabla^k$ to
$(2.1)_1,(2.1)_2$ and then multiplying the resulting equations by
$\nabla^kn_1,\nabla^ku_1$ respectively, summing up and integrating
over $\mathbb{R}^3$, one has
$$\arraycolsep=1.5pt  \begin{array}[b]{rl}
&\D\frac{1}{2}\frac{d}{dt}\int|\nabla^k(n_1,u_1)|^2+\|\nabla^ku_1\|_{L^2}^2-\int\nabla^ku_1\nabla^k\nabla\phi\\[3mm]
=&\D-\int\nabla^kn_1\nabla^k(u_1\cdot \nabla n_1+n_1{\rm
div}u_1)+\nabla^ku_1\nabla^k(u_1\cdot \nabla u_1+h(n_1)\nabla n_1)\\[3mm]
=&\D\!-\!\int\nabla^k(u_1\cdot\nabla
n_1)\nabla^kn_1\!-\!\nabla^k(u_1\cdot\nabla
u_1)\nabla^ku_1\!-\!\nabla^k(n_1{\rm
div}u_1)\nabla^kn_1\!-\!\nabla^k(h(n_1)\nabla n_1)\nabla^k
u_1\\[3mm]
:=&\D I_1+I_2+I_3+I_4.
\end{array}\eqno(2.5)
$$

We shall first estimate each term in the right hand side of (2.5).
By H\"{o}lder's inequalities and Lemma 2.1, we get
$$\arraycolsep=1.5pt  \begin{array}[b]{rl}
 I_1=&\D-\int\sum_{0\leq l\leq
k}C_k^l\nabla^{k-l}u_1\cdot\nabla\nabla^ln_1\nabla^kn_1\leq\sum_{0\leq
l\leq
k}\|\nabla^{k-l}u_1\nabla\nabla^ln_1\|_{L^{6/5}}\|\nabla^kn_1\|_{L^6}\\[3mm]
\leq &\D\sum_{0\leq l\leq
k}\|\nabla^{k-l}u_1\nabla\nabla^ln_1\|_{L^{6/5}}\|\nabla^{k+1}n_1\|_{L^2}.
\end{array}
\eqno(2.6)
$$

When $0\leq l\leq[\frac{k}{2}]$, by H\"{o}lder's inequality and
Lemma 2.1, we have
$$\arraycolsep=1.5pt  \begin{array}[b]{rl}
\|\nabla^{k-l}u_1\nabla\nabla^ln_1\|_{L^{6/5}}\leq &\D
\|\nabla^{k-l}u_1\|_{L^2}\|\nabla^{l+1}n_1\|_{L^3}\\[3mm]
\leq
&\D\|u_1\|_{L^2}^{\frac{l}{k}}\|\nabla^ku_1\|_{L^2}^{1-\frac{l}{k}}\|\nabla^\alpha
n_1\|^{1-\frac{l}{k}}\|\nabla^{k+1}n_1\|_{L^2}^{\frac{l}{k}}\\[3mm]
\leq &\D\delta(\|\nabla^{k+1}n_1\|_{L^2}+\|\nabla^ku_1\|_{L^2}),
\end{array}
\eqno(2.7)
$$
where $\alpha$ satisfies
$$
l+\frac{3}{2}=\alpha(1-\frac{l}{k})+(k+1)\frac{l}{k},
$$
which gives $\alpha=\frac{3k-2l}{2k-2l}\in[\frac{3}{2},3)$ since
$l\leq\frac{k}{2}$.

When $[\frac{k}{2}]+1\leq l\leq k$, by H\"{o}lder's inequality and
Lemma 2.1 again, we obtain
$$\arraycolsep=1.5pt  \begin{array}[b]{rl}
\|\nabla^{k-l}u_1\nabla\nabla^ln_1\|_{L^{6/5}}\leq &\D
\|\nabla^{k-l}u_1\|_{L^3}\|\nabla^{l+1}n_1\|_{L^2}\\[3mm]
\leq
&\D\|n_1\|_{L^2}^{\frac{k-l}{k+1}}\|\nabla^{k+1}n_1\|_{L^2}^{\frac{l+1}{k+1}}\|\nabla^\alpha
u_1\|^{1-\frac{l+1}{k+1}}\|\nabla^{k+1}u_1\|_{L^2}^{\frac{l-1}{k+1}}\\[3mm]
\leq &\D\delta(\|\nabla^{k+1}n_1\|_{L^2}+\|\nabla^ku_1\|_{L^2}),
\end{array}
\eqno(2.8)
$$
where $\alpha$ satisfies
$$
k-l+\frac{1}{2}=\alpha\frac{l+1}{k+1}+k\frac{k-l}{k+1},
$$
which implies $\alpha=\frac{3k-2l+1}{2l+2}\in[\frac{1}{2},3)$ since
$l\geq\frac{k+1}{2}$.

From (2.6), (2.7) and (2.8), one has
$$
I_1\leq
\delta(\|\nabla^{k+1}n_1\|_{L^2}+\|\nabla^ku_1\|_{L^2}).\eqno(2.9)
$$

For $I_2$, using Lemma 2.1 and H\"{o}lder's inequality, we get
$$\arraycolsep=1.5pt  \begin{array}[b]{rl}
I_2=&\D-\int([\nabla^k,u_1]\nabla
u_1+u_1\nabla\nabla^ku_1)\nabla^ku_1\leq\|\nabla
u_1\|_{L^\infty}\|\nabla^ku_1\|_{L^2}^2-\frac{1}{2}\int
u_1\nabla(\nabla^ku_1\nabla^ku_1)\\[3mm]
\leq &\D\|\nabla
u_1\|_{L^\infty}\|\nabla^ku_1\|_{L^2}^2+\frac{1}{2}\int{\rm
div}u_1\nabla^ku_1\cdot\nabla^ku_1\leq
\delta\|\nabla^ku_1\|_{L^2}^2.
\end{array}\eqno(2.10)
$$

For $I_3$,
$$\arraycolsep=1.5pt  \begin{array}[b]{rl}
I_3=&\D-\int\nabla^k(n_1{\rm div}u_1)\nabla^kn_1\\[3mm]
=&\D-\int\sum_{0\leq l\leq k-1}C_k^l\nabla^{k-l}n_1\nabla^l{\rm
div}u_1\nabla^kn_1-\int n_1{\rm div}\nabla^ku_1\nabla^kn_1\\[3mm]
:=&\D I_{31}+I_{32}.
\end{array}
\eqno(2.11)
$$

First, we estimate $I_{31}$. By H\"{o}lder's inequality, Lemma 2.1
and Cauchy's inequality, we obtain
$$\arraycolsep=1.5pt  \begin{array}[b]{rl}
I_{31}=&\D-\int\sum_{0\leq l\leq
k-1}C_k^l\nabla^{k-l}n_1\nabla^l{\rm div}u_1\nabla^kn_1\\[3mm]
\leq &\D C\sum_{0\leq l\leq k-1}\|\nabla^{k-l}n_1\nabla^l{\rm
div}u_1\|_{L^{6/5}}\|\nabla^{k+1}n_1\|_{L^2}.
\end{array}
\eqno(2.12)
$$

When $0\leq l\leq[\frac{k}{2}]$, using Lemma 2.1 and H\"{o}lder's
inequality, we have
$$\arraycolsep=1.5pt  \begin{array}[b]{rl}
\|\nabla^{k-l}n_1\nabla^l{\rm div}u_1\|_{L^{6/5}}\leq &\D
C\|\nabla^{k-l}n_1\|_{L^2}\|\nabla^{l+1}u_1\|_{L^3}\\[3mm]
\leq &\D
C\|n_1\|_{L^2}^{\frac{l+1}{k+1}}\|\nabla^{k+1}n_1\|_{L^2}^{\frac{k-l}{k+1}}\|\nabla^\alpha
u_1\|_{L^2}^{\frac{k-l}{k+1}}\|\nabla^ku_1\|_{L^2}^{\frac{l+1}{k+1}}\\[3mm]
\leq &\D C\delta(\|\nabla^{k+1}n_1\|_{L^2}+\|\nabla^ku_1\|_{L^2},
\end{array}
\eqno(2.13)
$$
where $\alpha$ satisfies
$$
l+\frac{3}{2}=\alpha\frac{k-l}{k+1}+k\frac{l+1}{k+1},
$$
which yields $\alpha=\frac{k+2l+3}{2k-2l}\in(\frac{1}{2},3)$ since
$l\leq\frac{k}{2}$.

When $[\frac{k}{2}]+1\leq l\leq k-1$, using Lemma 2.1 and
H\"{o}lder's inequality, we have
$$\arraycolsep=1.5pt  \begin{array}[b]{rl}
\|\nabla^{k-l}n_1\nabla^l{\rm div}u_1\|_{L^{6/5}}\leq &\D
C\|\nabla^{k-l}n_1\|_{L^3}\|\nabla^{l+1}u_1\|_{L^2}\\[3mm]
\leq &\D C\|\nabla^\alpha
n_1\|_{L^2}^{\frac{l+1}{k}}\|\nabla^{l+1}n_1\|_{L^2}^{\frac{k-1-l}{k}}\|u_1\|_{L^2}^{\frac{k-1-l}{k}}
\|\nabla^ku_1\|_{L^2}^{\frac{l+1}{k}}\\[3mm]
\leq &\D C\delta(\|\nabla^{k+1}n_1\|_{L^2}+\|\nabla^ku_1\|_{L^2}),
\end{array}
\eqno(2.14)
$$
where $\alpha$ satisfies
$$
k-l+\frac{1}{2}=\alpha\frac{l+l}{k}+(k+1)\frac{k-l-1}{k},
$$
which yields $\alpha=1+\frac{k}{2l+2}\in(\frac{3}{2},3)$ since
$l\geq\frac{k+1}{2}$.

From (2.12), (2.13) and (2.14), we get
$$
I_{31}\leq
C\delta(\|\nabla^{k+1}n_1\|_{L^2}^2+\|\nabla^ku_1\|_{L^2}^2).\eqno(2.15)
$$

For $I_{32}$, By H\"{o}lder's inequality, Lemma 2.1 and Cauchy's
inequality, we obtain
$$\arraycolsep=1.5pt  \begin{array}[b]{rl}
I_{32}=&\D-\int n_1{\rm div}\nabla^ku_1\nabla^kn_1\\[3mm]
=&\D-\int n_1{\rm div}(\nabla^ku_1\nabla^kn_1)+\int
n_1\nabla^{k+1}n_1\nabla^ku_1\\[3mm]
\leq &\D C\|\nabla
n_1\|_{L^3}\|\nabla^ku_1\|_{L^2}\|\nabla^kn_1\|_{L^6}+\|n_1\|_{L^\infty}\|\nabla^{k+1}n_1\|_{L^2}\|\nabla^ku_1\|_{L^2}\\[3mm]
\leq &\D
C\delta(\|\nabla^{k+1}n_1\|_{L^2}^2+\|\nabla^ku_1\|_{L^2}^2).
\end{array}
\eqno(2.16)
$$

Thus, (2.11), (2.15) and (2.16) imply
$$
I_3\leq
C\delta(\|\nabla^{k+1}n_1\|_{L^2}^2+\|\nabla^ku_1\|_{L^2}^2).\eqno(2.17)
$$

Next, we will estimate $I_4$.
$$
\arraycolsep=1.5pt  \begin{array}[b]{rl}
I_4=&\D-\int\nabla^k(h(n_1)\nabla n_1)\nabla^ku_1\\[3mm]
=&\D-\int\sum_{0\leq l\leq
k}C_k^l\nabla^{k-l}h(n_1)\nabla^{l+1}n_1\nabla^ku_1
+h(n_1)\nabla^{k+1}n_1\cdot\nabla^ku_1\\[3mm]
:=&\D I_{41}+I_{42}.
\end{array}
\eqno(2.18)
$$

For $I_{41}$, By H\"{o}lder's inequality and Lemma 2.1, we obtain
$$\arraycolsep=1.5pt  \begin{array}[b]{rl}
I_{41}=&\D-\int\sum_{0\leq l\leq
k}C_k^l\nabla^{k-l}h(n_1)\nabla^{l+1}n_1\nabla^ku_1\\[3mm]
\leq &\D
C\|\nabla^{k-l}n_1\nabla^{l+1}n_1\|_{L^2}\|\nabla^ku_1\|_{L^2}.
\end{array}
\eqno(2.19)
$$

When $0\leq l\leq[\frac{k}{2}]$, by using H\"{o}lder's inequality
and Lemma 2.1, we get
$$\arraycolsep=1.5pt  \begin{array}[b]{rl}
&\D\|\nabla^{k-l}h(n_1)\nabla^{l+1}n_1\|_{L^2}
\leq \|\nabla^{k-l}h(n_1)\|_{L^6}\|\nabla^{l+1}n_1\|_{L^3}\\[3mm]
 \leq
&\D
C\|\nabla^{k-l}h(n_1)\|_{L^2}^{\frac{l}{k+1}}\|\nabla^{k+1}h(n_1)\|_{L^2}^{1-\frac{l}{k+1}}\|\nabla^\alpha
n_1\|_{L^2}^{1-\frac{l}{k+1}}\|\nabla^{k+1}n_1\|_{L^2}^{\frac{l}{k+1}}\\[3mm]
\leq &\D
C\|\nabla^{k-l}n_1\|_{L^2}^{\frac{l}{k+1}}\|\nabla^{k+1}n_1\|_{L^2}^{1-\frac{l}{k+1}}\|\nabla^\alpha
n_1\|_{L^2}^{1-\frac{l}{k+1}}\|\nabla^{k+1}n_1\|_{L^2}^{\frac{l}{k+1}}\\[3mm]
\leq &\D C\delta\|\nabla^{k+1}n_1\|_{L^2},
\end{array}
\eqno(2.20)
$$
where $\alpha$ satisfies
$$
l+\frac{3}{2}=\alpha(1-\frac{l}{k+1})+l,
$$
which implies $\alpha=\frac{3k+3}{2k-2l+2}\in[\frac{3}{2},3)$, since
$l\leq\frac{k}{2}$.

When $[\frac{k}{2}]+1\leq l\leq k-1$, by H\"{o}lder's inequality and
Lemma 2.1, we get
$$\arraycolsep=1.5pt  \begin{array}[b]{rl}
&\D\|\nabla^{k-l}h(n_1)\nabla^{l+1}n_1\|_{L^2}
\leq \|\nabla^{k-l}h(n_1)\|_{L^3}\|\nabla^{l+1}n_1\|_{L^6}\\[3mm]
 \leq
&\D C\|\nabla^\alpha
h(n_1)\|_{L^2}^{\frac{l}{k-1}}\|\nabla^{k+1}h(n_1)\|_{L^2}^{1-\frac{l}{k-1}}\|\nabla^2
n_1\|_{L^2}^{1-\frac{l}{k-1}}\|\nabla^{k+1}n_1\|_{L^2}^{\frac{l}{k-1}}\\[3mm]
\leq &\D C\|\nabla^\alpha
n_1\|_{L^2}^{\frac{l}{k-1}}\|\nabla^{k+1}n_1\|_{L^2}^{1-\frac{l}{k-1}}\|\nabla^2
n_1\|_{L^2}^{1-\frac{l}{k-1}}\|\nabla^{k+1}n_1\|_{L^2}^{\frac{l}{k-1}}\\[3mm]
\leq &\D C\delta\|\nabla^{k+1}n_1\|_{L^2},
\end{array}
\eqno(2.21)
$$
where $\alpha$ satisfies
$$
k-l+\frac{1}{2}=\alpha\frac{l}{k-1}+(k+1)(1-\frac{l}{k-1}),
$$
which implies $\alpha=2+\frac{-k+1}{2l}\in[\frac{3}{2},3)$ since
$l\geq\frac{k+1}{2}$.

Thus, from (2.18), (2.19), (2.20) and (2.21), we deduce that
$$
I_4\leq
C\delta(\|\nabla^{k+1}n_1\|_{L^2}^2+\|\nabla^ku_1\|_{L^2}^2).\eqno(2.22)
$$

Hence, for $n_1$ and $u_1$, we have
$$\arraycolsep=1.5pt  \begin{array}[b]{rl}
&\D\frac{1}{2}\frac{d}{dt}\int|\nabla^k(n_1,u_1)|^2+\|\nabla^ku_1\|_{L^2}^2-\int\nabla^ku_1\nabla^k\nabla\phi\\[3mm]
\leq &\D
C\delta(\|\nabla^{k+1}n_1\|_{L^2}^2+\|\nabla^ku_1\|_{L^2}^2).
\end{array}\eqno(2.23)
$$

In the same way, we can get the following estimates for $n_2$ and
$u_2$, that is,
$$
\arraycolsep=1.5pt  \begin{array}[b]{rl}
&\D\frac{1}{2}\frac{d}{dt}\int|\nabla^k(n_2,u_2)|^2+\|\nabla^ku_2\|_{L^2}^2+\int\nabla^ku_2\nabla^k\nabla\phi\\[3mm]
\leq &\D
C\delta(\|\nabla^{k+1}n_2\|_{L^2}^2+\|\nabla^ku_2\|_{L^2}^2).
\end{array}\eqno(2.24)
$$

Finally, we will turn to estimate the last term in left hand side of
(2.23) and (2.24). Since $n_1$ and $n_2$ is coupled in Poisson
equation, we will estimate them  simultaneously as follows.
$$\arraycolsep=1.5pt  \begin{array}[b]{rl}
&\D-\int\nabla^k\nabla\phi\cdot\nabla^ku_1+\int\nabla^k\nabla\phi\cdot\nabla^ku_2\\[3mm]
=&\D\int\nabla^k({\rm div}u_1)\nabla^k\phi-\int\nabla^k({\rm
div}u_2)\nabla^k\phi\\[3mm]
=&\D-\int\nabla^k[\partial_tn_1+{\rm
div}(n_1u_1)\nabla^k\phi+\int\nabla^k[\partial_tn_2+{\rm
div}(n_2u_2)\nabla^k\phi\\[3mm]
=&\D-\int\nabla^k\partial_t(n_1-n_2)\nabla^k\phi-\int\nabla^k({\rm
div}(n_1u_1))\nabla^k\phi+\int\nabla^k({\rm
div}(n_2u_2))\nabla^k\phi\\[3mm]
=&\D-\int\nabla^k\partial_t\Delta\phi\nabla^k\phi-\int\nabla^k({\rm
div}(n_1u_1))\nabla^k\phi+\int\nabla^k({\rm
div}(n_2u_2))\nabla^k\phi\\[3mm]
=&\D\frac{1}{2}\frac{d}{dt}\|\nabla^k\nabla\phi\|_{L^2}^2+\int\nabla^k(n_1u_1)\nabla^k\nabla\phi-\int\nabla^k(n_2u_2)\nabla^k\nabla\phi\\[3mm]
:=&\D\frac{1}{2}\frac{d}{dt}\|\nabla^k\nabla\phi\|_{L^2}^2+I_{51}+I_{52}
\end{array}
\eqno(2.25)
$$

Now we will estimate $I_{51}$ and $I_{52}$.

When $k=0$, by H\"{o}lder's inequality, Sobolev's inequality and
Cauchy's inequality, we have
$$\arraycolsep=1.5pt  \begin{array}[b]{rl}
\D\int n_1u_1\nabla\phi\leq &\D
C\|\nabla\phi\|_{L^6}\|u_1\|_{L^2}\|n_1\|_{L^3}
\leq C\|\nabla\nabla\phi\|_{L^2}\|u_1\|_{L^2}\|n_1\|_{L^3}\\[3mm]
\leq &\D C\delta(\|\nabla\nabla\phi\|_{L^2}+\|u\|_{L^2}^2).
\end{array}
\eqno(2.26)
$$

Similarly, for $k=1$, we get
$$\arraycolsep=1.5pt  \begin{array}[b]{rl}
\D\int\nabla(n_1u_1)\nabla(\nabla\phi)=&\D-\int(n_1u_1)\nabla^2\nabla\phi\leq
C\|\nabla^2\nabla\phi\|_{L^2}\|u\|_{L^6}\|n\|_{L^3}\\[3mm]
\leq &\D C\|\nabla^2\nabla\phi\|_{L^2}\|\nabla
u_1\|_{L^2}\|n\|_{L^3}\\[3mm]
\leq &\D C\delta(\|\nabla^2\nabla\phi\|_{L^2}+\|\nabla
u_1\|_{L^2}^2);
\end{array}
\eqno(2.27)
$$
and for $k=2$, we have
$$\arraycolsep=1.5pt  \begin{array}[b]{rl}
\D\int\nabla^2(n_1u_1)\nabla^2(\nabla\phi)=&\D-\int\nabla(n_1u_1)\nabla^3\nabla\phi
\leq \|\nabla^3\nabla\phi\|_{L^2}\|\sum_{0\leq
l\leq1}\nabla^{1-l}n_1\nabla^lu_1\|_{L^2}\\[3mm]
\leq &\D C\|\nabla^3\nabla\phi\|_{L^2}\|\nabla^\alpha
n_1\|_{L^2}^{\frac{l+1}{2}}\|\nabla^3n_1\|_{L^2}^{1-\frac{l+1}{2}}\|u_1\|_{L^2}^{1-\frac{l+1}{2}}\|\nabla^2u_1\|_{L^2}^{\frac{l+1}{2}}\\[3mm]
\leq &\D
C\delta(\|\nabla^3\nabla\phi\|_{L^2}^2+\|\nabla^3n_1\|_{L^2}^2+\|\nabla^2u_1\|_{L^2}^2),
\end{array}
\eqno(2.28)
$$
where
$$
\alpha=\frac{l}{l+1},\ l=0,1.
$$

In the same way, one can obtain the estimate of $I_{52}$. Hence,
from (2.25) to (2.28), we have
$$\arraycolsep=1.5pt  \begin{array}[b]{rl}
I_{51}\!+\!I_{52}\!\geq\!
-C\delta(\|\nabla^{k+1}n_1\|_{L^2}^2\!+\!\|\nabla^ku_1\|_{L^2}^2\!+\!\|\nabla^{k+1}n_2\|_{L^2}^2\!+\!\|\nabla^ku_2\|_{L^2}^2
\!+\!\|\nabla^{k+1}\nabla\phi\|_{L^2}^2).
\end{array}
\eqno(2.29)
$$

Combining (2.23), (2.24), (2.25) and (2.29), we deduce that
$$\arraycolsep=1.5pt  \begin{array}[b]{rl}
&\D\frac{1}{2}\frac{d}{dt}\int|\nabla^k(n_1,u_1,n_2,u_2,\nabla\phi)|^2+\|\nabla^k(u_1,u_2)\|_{L^2}^2\\[3mm]
\leq &\D
C\delta(\|\nabla^{k+1}n_1\|_{L^2}^2+\|\nabla^ku_1\|_{L^2}^2+\|\nabla^{k+1}n_2\|_{L^2}^2
+\|\nabla^ku_2\|_{L^2}^2 +\|\nabla^{k+1}\nabla\phi\|_{L^2}^2).
\end{array}
\eqno(2.30)
$$

This proves Lemma 2.10. \qed

Next, we derive the second type of energy estimates excluding
$n_1,u_1$ and $n_2,u_2$ themselves.

\bigbreak \noindent\textbf{Lemma 2.11.}  Assume that $0\leq k\leq2$,
then we have
$$\arraycolsep=1.5pt  \begin{array}[b]{rl}
&\D\frac{1}{2}\frac{d}{dt}\int|\nabla^{k+1}(n_1,u_1,n_2,u_2,\nabla\phi)|^2+\|\nabla^{k+1}(u_1,u_2)\|_{L^2}^2\\[3mm]
\leq &\D
C\delta(\|\nabla^{k+1}n_1\|_{L^2}^2+\|\nabla^{k+1}u_1\|_{L^2}^2+\|\nabla^{k+1}n_2\|_{L^2}^2
+\|\nabla^{k+1}u_2\|_{L^2}^2 +\|\nabla^{k+1}\nabla\phi\|_{L^2}^2).
\end{array}
\eqno(2.31)
$$
\noindent\textbf{Proof.} For $0\leq k\leq 2$, applying
$\nabla^{k+1}$ to $(2.1)_1,(2.1)_2$ and then multiplying the
resulting equations by $\nabla^{k+1}n_1,\nabla^{k+1}u_1$
respectively, summing up and integrating over $\mathbb{R}^3$, one
has
$$\arraycolsep=1.5pt  \begin{array}[b]{rl}
&\D\frac{1}{2}\frac{d}{dt}\int|\nabla^{k+1}(n_1,u_1)|^2+\|\nabla^{k+1}u_1\|_{L^2}^2-\int\nabla^{k+1}u_1\cdot\nabla^{k+1}\nabla\phi\\[3mm]
=&\D-\int\nabla^{k+1}n_1\nabla^{k+1}(u_1\cdot \nabla n_1+n_1{\rm
div}u_1)+\nabla^{k+1}u_1\nabla^{k+1}(u_1\cdot \nabla u_1+h(n_1)\nabla n_1)\\[3mm]
=&\D\!-\!\int[\nabla^{k+1}(u_1\cdot\nabla
n_1)\nabla^{k+1}n_1\!+\!\nabla^{k+1}(u_1\cdot\nabla
u_1)\nabla^{k+1}u_1\!]\\[3mm]
&\D-\int[\nabla^{k+1}(n_1{\rm
div}u_1)\nabla^{k+1}n_1\!+\!\nabla^{k+1}(h(n_1)\nabla
n_1)\nabla^{k+1}
u_1]\\[3mm]
:=&\D J_1+J_2.
\end{array}\eqno(2.32)
$$

Now we shall estimate $J_1$ and $J_2$. By Lemma 2.2, H\"{o}lder's
inequality and Cauchy's inequality, we get
$$
\arraycolsep=1.5pt  \begin{array}[b]{rl}
J_1=&\D-\int \nabla^{k+1}(u_1\cdot\nabla
n_1)\nabla^{k+1}n_1\!+\!\nabla^{k+1}(u_1\cdot\nabla
u_1)\nabla^{k+1}u_1\\[3mm]
=&\D-\int[\nabla^{k+1},u_1]\cdot\nabla
n_1\nabla^{k+1}n_1+([\nabla^{k+1},u_1],\nabla
u_1)\cdot\nabla^{k+1}u_1\\[3mm]
&\D-\int
u_1\cdot\nabla\nabla^{k+1}n_1\nabla^{k+1}n_1+(u_1\cdot\nabla\nabla^{k+1}u_1)\cdot\nabla^{k+1}u_1\\[3mm]
\leq &\D C(\|\nabla
u_1\|_{L^\infty}\|\nabla^{k+1}n_1\|_{L^2}+\|\nabla^{k+1}u_1\|_{L^2}\|\nabla
n_1\|_{L^\infty})\|\nabla^{k+1}n_1\|_{L^2}\\[3mm]
&\D+\|\nabla
u_1\|_{L^\infty}\|\nabla^{k+1}u_1\|_{L^2}-\frac{1}{2}\int
u_1\cdot\nabla(\nabla^{k+1}n_1\nabla^{k+1}n_1+\nabla^{k+1}u_1\cdot\nabla^{k+1}u_1)\\[3mm]
\leq &\D
C\|\nabla(n_1,u_1)\|_{L^\infty}\|\nabla^{k+1}(n_1,u_1)\|_{L^2}^2+\frac{1}{2}{\rm
div}u_1\nabla^{k+1}n_1\nabla^{k+1}n_1+{\rm
div}u_1\nabla^{k+1}u_1\cdot\nabla^{k+1}u_1\\[3mm]
\leq &\D
C\delta(\|\nabla^{k+1}n_1\|_{L^2}^2+\|\nabla^{k+1}u_1\|_{L^2}^2).
\end{array}
\eqno(2.33)
$$

In the same way, one can deduce that
$$
J_2\leq
C\delta(\|\nabla^{k+1}n_1\|_{L^2}^2+\|\nabla^{k+1}u_1\|_{L^2}^2).\eqno(2.34)
$$

Thus we have
$$\arraycolsep=1.5pt  \begin{array}[b]{rl}
&\D\frac{1}{2}\frac{d}{dt}\int|\nabla^{k+1}(n_1,u_1)|^2+\|\nabla^{k+1}u_1\|_{L^2}^2-\int\nabla^{k+1}u_1\nabla^{k+1}\nabla\phi\\[3mm]
\leq &\D
C\delta(\|\nabla^{k+1}n_1\|_{L^2}^2+\|\nabla^{k+1}u_1\|_{L^2}^2).
\end{array}\eqno(2.35)
$$

The similar estimate of $n_2,u_2$ is
$$\arraycolsep=1.5pt  \begin{array}[b]{rl}
&\D\frac{1}{2}\frac{d}{dt}\int|\nabla^{k+1}(n_2,u_2)|^2+\|\nabla^{k+1}u_2\|_{L^2}^2+\int\nabla^{k+1}u_2\nabla^{k+1}\nabla\phi\\[3mm]
\leq &\D
C\delta(\|\nabla^{k+1}n_2\|_{L^2}^2+\|\nabla^{k+1}u_2\|_{L^2}^2).
\end{array}\eqno(2.36)
$$

Finally, we give the estimates of the last terms in the left hand
side of (2.35) and (2.36) as follows.
$$
\arraycolsep=1.5pt  \begin{array}[b]{rl}
&\D-\int\nabla^{k+1}\nabla\phi\cdot\nabla^{k+1}u_1+\int\nabla^{k+1}\nabla\phi\cdot\nabla^{k+1}u_2\\[3mm]
=&\D\int\nabla^{k+1}({\rm
div}u_1)\nabla^{k+1}\phi-\int\nabla^{k+1}({\rm
div}u_2)\nabla^{k+1}\phi\\[3mm]
=&\D-\int\nabla^{k+1}[\partial_tn_1+{\rm
div}(n_1u_1)]\nabla^{k+1}\phi+\int\nabla^{k+1}[\partial_tn_2+{\rm
div}(n_2u_2)]\nabla^{k+1}\phi\\[3mm]
=&\D-\int\nabla^{k+1}\partial_t(n_1-n_2)\nabla^{k+1}\phi-\int\nabla^{k+1}({\rm
div}(n_1u_1))\nabla^{k+1}\phi+\int\nabla^{k+1}({\rm
div}(n_2u_2))\nabla^{k+1}\phi\\[3mm]
=&\D-\int\nabla^{k+1}\partial_t\Delta\phi\nabla^{k+1}\phi-\int\nabla^{k+1}({\rm
div}(n_1u_1))\nabla^{k+1}\phi+\int\nabla^{k+1}({\rm
div}(n_2u_2))\nabla^{k+1}\phi\\[3mm]
=&\D\frac{1}{2}\frac{d}{dt}\|\nabla^{k+1}\nabla\phi\|_{L^2}^2+\int\nabla^{k+1}(n_1u_1)\nabla^{k+1}\nabla\phi
-\int\nabla^{k+1}(n_2u_2)\nabla^{k+1}\nabla\phi\\[3mm]
:=&\D\frac{1}{2}\frac{d}{dt}\|\nabla^{k+1}\nabla\phi\|_{L^2}^2+J_3+J_4.
\end{array}
\eqno(2.37)
$$

Using H\"{o}lder's inequality and Lemma 2.2 and Cauchy's inequality,
we obtain
$$\arraycolsep=1.5pt  \begin{array}[b]{rl}
J_3=&\D\int\nabla^{k+1}(n_1u_1)\cdot\nabla^{k+1}\nabla\phi \leq
C\|\nabla^{k+1}\nabla\phi\|_{L^2}\|\nabla^{k+1}(n_1u_1)\|_{L^2}\\[3mm]
\leq &\D C\|\nabla^{k+1}\nabla\phi\|_{L^2}(\|n_1\|_{L^\infty}\|\nabla^{k+1}u_1\|_{L^2}
+\|u_1\|_{L^\infty}\|\nabla^{k+1}n_1\|_{L^2})\\[3mm]
\leq &\D
C\delta(\|\nabla^{k+1}u_1\|_{L^2}^2+\|\nabla^{k+1}n_1\|_{L^2}^2+\|\nabla^{k+1}\nabla\phi\|_{L^2}^2).
\end{array}
\eqno(2.38)
$$

Similarly, we have
$$\arraycolsep=1.5pt  \begin{array}[b]{rl}
J_4=&\D\int\nabla^{k+1}(n_2u_2)\cdot\nabla^{k+1}\nabla\phi \leq
C\|\nabla^{k+1}\nabla\phi\|_{L^2}\|\nabla^{k+1}(n_2u_2)\|_{L^2}\\[3mm]
\leq &\D
C\delta(\|\nabla^{k+1}u_2\|_{L^2}^2+\|\nabla^{k+1}n_2\|_{L^2}^2+\|\nabla^{k+1}\nabla\phi\|_{L^2}^2).
\end{array}
\eqno(2.39)
$$

Hence, plugging (2.33), (2.34), (2.37), (2.38) and (2.39) into
(2.32), we deduce that (2.31). This proves Lemma 2.11. \qed

Now, we shall recover the dissipation estimate for $n_1,n_2$.

\bigbreak \noindent\textbf{Lemma 2.12.} Assume that $0\leq k\leq 2$,
then we have
$$
\arraycolsep=1.5pt  \begin{array}[b]{rl}
&\D\frac{d}{dt}\left\{\int\nabla^ku_1\cdot\nabla\nabla^kn_1+\nabla^ku_2\cdot\nabla\nabla^kn_2\right\}
+C\|\nabla^{k+1}(n_1,n_2,\nabla\phi)\|_{L^2}^2\\[3mm]
\leq &\D
C(\|\nabla^ku_1\|_{L^2}^2+\|\nabla^{k+1}u_1\|_{L^2}^2+\|\nabla^ku_2\|_{L^2}^2+\|\nabla^{k+1}u_2\|_{L^2}^2).
\end{array}
\eqno(2.40)
$$
\noindent\textbf{Proof.} Let $0\leq k\leq 2$. Applying $\nabla^k$ to
$(2.1)_2$ and then multiplying the resulting equality by
$\nabla\nabla^kn_1$, we have
$$\arraycolsep=1.5pt  \begin{array}[b]{rl}
\|\nabla^{k+1}n_1\|_{L^2}^2-\int\nabla\nabla^kn_1\nabla^k\nabla\phi
\leq &\D
-\int\nabla^k\partial_tu_1\cdot\nabla\nabla^kn_1+C\|\nabla^ku_1\|_{L^2}\|\nabla^{k+1}n_1\|_{L^2}\\[3mm]
&\D+\|\nabla^k(u_1\cdot\nabla u_1+h(n_1)\nabla
n_1)\|_{L^2}\|\nabla^{k+1}n_1\|_{L^2}.
\end{array}
\eqno(2.41)
$$

First, we estimate the first term in the right hand side of (2.39).
$$\arraycolsep=1.5pt  \begin{array}[b]{rl}
&\D-\int\nabla^ku_1\partial_tu_1\cdot\nabla\nabla^kn_1\\[3mm]
=&\D-\frac{d}{dt}\int\nabla^k\cdot\nabla\nabla^kn_1-\int\nabla^k{\rm
div}u_1\cdot\nabla^k\partial_tn_1\\[3mm]
=&\D-\frac{d}{dt}\int\nabla^ku_1\cdot\nabla\nabla^kn_1+\|\nabla^k{\rm
div}u_1\|_{L^2}^2+\int\nabla^ku_1{\rm
div}\cdot\nabla^k(u_1\cdot\nabla n_1+n_1{\rm div}u_1),
\end{array}
\eqno(2.42)
$$

Next, we shall estimate the last two terms in (2.40) by
$$
\arraycolsep=1.5pt  \begin{array}[b]{rl} \D\int\nabla^k{\rm
div}u_1\cdot\nabla^k(u_1\cdot\nabla n_1)=&\D\int\sum_{0\leq l\leq
k}C_k^l\nabla^lu_1\cdot\nabla\nabla^{k-l}n_1\cdot\nabla^k{\rm
div}u_1\\[3mm]
\leq &\D C\sum_{0\leq l\leq
k}\|\nabla^lu_1\cdot\nabla\nabla^{k-l}n_1\|_{L^2}\|\nabla^{k+1}u_1\|_{L^2}.
\end{array}
\eqno(2.43)
$$

If $l=0$, then
$$\arraycolsep=1.5pt  \begin{array}[b]{rl}
\|u_1\cdot\nabla\nabla^kn_1\|_{L^2}\|\nabla^{k+1}u_1\|_{L^2}\leq &\D
C\|u_1\|_{L^\infty}\|\nabla^{k+1}n_1\|_{L^2}\|\nabla^{k+1}u_1\|_{L^2}\\[3mm]
\leq &\D
C\delta(\|\nabla^{k+1}n_1\|_{L^2}^2+\|\nabla^{k+1}u_1\|_{L^2}^2).
\end{array}
\eqno(2.44)
$$

If $1\leq l\leq[k/2]$, using H\"{o}lder's inequality and Lemma 2.1,
we get
$$\arraycolsep=1.5pt  \begin{array}[b]{rl}
\|\nabla^lu_1\cdot\nabla\nabla^{k-l}n_1\|_{L^2}\leq
&\D C\|\nabla^{k+1-l}\|_{L^6}\|\nabla^lu_1\|_{L^3}\\[3mm]
\leq &\D
C\|n_1\|_{L^2}^{\frac{l-1}{k+1}}\|\nabla^{k+1}n_1\|_{L^2}^{\frac{k-l+2}{k+1}}\|\nabla^\alpha
u_1\|_{L^2}^{\frac{k-l+2}{k+1}}\|\nabla^{k+1}u_1\|_{L^2}^{\frac{l-1}{k+1}}\\[3mm]
\leq &\D
C\delta(\|\nabla^{k+1}n_1\|_{L^2}+\|\nabla^{k+1}u_1\|_{L^2}),
\end{array}
\eqno(2.45)
$$
where $\alpha=\frac{3k+3}{2k-2l+4}\in[3/2,3),\ {\rm since}\ l\leq
k/2$.

If $[k/2]+1\leq l\leq k$, using H\"{o}lder's inequality and Lemma
2.1 again, we obtain
$$\arraycolsep=1.5pt  \begin{array}[b]{rl}
\|\nabla^lu_1\cdot\nabla\nabla^{k-l}n_1\|_{L^2}\leq
&\D C\|\nabla^{k+1-l}\|_{L^3}\|\nabla^lu_1\|_{L^63}\\[3mm]
\leq &\D C\|\nabla^\alpha
n_1\|_{L^2}^{\frac{l+1}{k+1}}\|\nabla^{k+1}n_1\|_{L^2}^{\frac{k-l}{k+1}}
\|u_1\|_{L^2}^{\frac{k-l}{k+1}}\|\nabla^{k+1}u_1\|_{L^2}^{\frac{l+1}{k+1}}\\[3mm]
\leq &\D
C\delta(\|\nabla^{k+1}n_1\|_{L^2}+\|\nabla^{k+1}u_1\|_{L^2}),
\end{array}
\eqno(2.46)
$$
where $\alpha=\frac{3k+3}{2l+2}\in[3/2,3),\ {\rm since}\ l\geq
\frac{k+1}{2}$.

Thus, from (2.44), (2.45) and (2.46), we obtain
$$
\int\nabla^k{\rm div}u_1\cdot\nabla^k(u_1\cdot\nabla n_1)\leq
C\delta(\|\nabla^{k+1}n_1\|_{L^2}+\|\nabla^{k+1}u_1\|_{L^2}).\eqno(2.47)
$$

Similarly, we also get
$$
\int\nabla^k{\rm div}u_1\cdot\nabla^k(n_1{\rm div}u_1)\leq
C\delta(\|\nabla^{k+1}n_1\|_{L^2}+\|\nabla^{k+1}u_1\|_{L^2}),\eqno(2.48)
$$
and
$$
\|\nabla^k(u_1\cdot\nabla u_1+h(n_1)\nabla n_1)\|_{L^2}\leq
C\delta(\|\nabla^{k+1}n_1\|_{L^2}+\|\nabla^{k+1}u_1\|_{L^2}).\eqno(2.49)
$$

Hence, by (2.40)-(2.49), we have
$$
\arraycolsep=1.5pt  \begin{array}[b]{rl}
&\D\frac{d}{dt}\int\nabla^ku_1\cdot\nabla\nabla^kn_1+C\|\nabla^{k+1}n_1\|_{L^2}-\int\nabla\nabla^kn_1\nabla^k\nabla\phi\\[3mm]
\leq &\D C(\|\nabla^ku_1\|_{L^2}^2+\|\nabla^{k+1}u_1\|_{L^2}^2).
\end{array}
\eqno(2.50)
$$

On the other hand, by a method similar to the above, we have
$$
\arraycolsep=1.5pt  \begin{array}[b]{rl}
&\D\frac{d}{dt}\int\nabla^ku_2\cdot\nabla\nabla^kn_2+C\|\nabla^{k+1}n_2\|_{L^2}+\int\nabla\nabla^kn_2\nabla^k\nabla\phi\\[3mm]
\leq &\D C(\|\nabla^ku_2\|_{L^2}^2+\|\nabla^{k+1}u_2\|_{L^2}^2).
\end{array}
\eqno(2.51)
$$

Finally, using the Poisson equation in (2.1), the second terms on
the left hand side of (2.50) and (2.51) can be estimated as
$$
-\int\nabla\nabla^kn_1\nabla^k\nabla\phi+\int\nabla\nabla^kn_2\nabla^k\nabla\phi
=\frac{1}{2}\|\nabla^{k+1}\nabla\phi\|_{L^2}^2. \eqno(2.52)
$$

Summing (2.50) and (2.51), and using (2.52), one has
$$
\arraycolsep=1.5pt  \begin{array}[b]{rl}
&\D\frac{d}{dt}\left\{\int\nabla^ku_2\cdot\nabla\nabla^kn_2+\nabla^ku_1\cdot\nabla\nabla^kn_1\right\}
+C\|\nabla^{k+1}(n_1,n_2,\nabla\phi)\|_{L^2}\\[3mm]
\leq &\D
C(\|\nabla^ku_1\|_{L^2}^2+\|\nabla^{k+1}u_1\|_{L^2}^2+\|\nabla^ku_2\|_{L^2}^2+\|\nabla^{k+1}u_2\|_{L^2}^2).
\end{array}
\eqno(2.53)
$$

This proves (2.40). \qed

\subsection*{ 2.3.\ Estimates in $\dot{H}^{-s}(\mathbb{R}^3)$}

\quad\quad The following lemma plays a key role in the proof of
Theorem 1.2. It shows an energy estimate of the solutions in the
negative Sobolev space $\dot{H}^{-s}(\mathbb{R}^3)$. Namely, we have

 \noindent\textbf{Lemma 2.13.} If
$\|n_{i0},u_{i0},\nabla\phi_0\|_{H^3}\ll1$ with $i=1,2$, for
$s\in(0,\frac{1}{2}]$, we have
$$
\frac{\rm d}{{\rm d}t}\|(n_i,u_i,\nabla\phi)\|_{\dot{H}^{-s}}
 \leq
 C(\|\nabla n_i\|_{H^1}^2+\|u_i\|_{H^2}^2)\|(n_i,u_i,\nabla\phi)\|_{\dot{H}^{-s}},\
 i=1,2; \eqno(2.54)
$$
and for $s\in(\frac{1}{2},\frac{3}{2})$, we have
$$\arraycolsep=1.5pt  \begin{array}[b]{rl}
 &\D \frac{\rm d}{{\rm d}t}\|(n_i,u_i,\nabla\phi)\|_{\dot{H}^{-s}}\\[3mm]
 \leq &\D
 C\left\{\|(n_i,u_i)|_{L^2}^{s-\frac{1}{2}}\|\nabla(n_i,u_i)\|_{H^1}^{\frac{5}{2}-s}
 +\|u_i\|_{L^2}\|n_i\|_{L^2}^{s-\frac{1}{2}}\|\nabla n_i\|_{L^2}^{\frac{3}{2}-s}\right\}\|(n_i,u_i,\nabla\phi)\|_{\dot{H}^{-s}},\
 i=1,2.
 \end{array}
 \eqno(2.55)
$$
\noindent\textbf{Proof.} Applying $\Lambda^{-s}$ to $(2.2)_1$,
$(2.2)_2$ and multiplying the resulting identity by
$\Lambda^{-s}n_1$, $\Lambda^{-s}u_1$, respectively, and integrating
over $\mathbb{R}^3$ by parts, we get
$$\arraycolsep=1.5pt  \begin{array}[b]{rl}
 &\D \frac{\rm d}{{\rm d}t}\int\left(|\Lambda^{-s}n_i|^2+|\Lambda^{-s}u_i|^2\right)
 +
 \int|\nabla\Lambda^{-s}u_i|^2+(-1)^i\int\Lambda^{-s}\nabla\phi\cdot\Lambda^{-s}u_i
\\[3mm]
 = &\D \int\Lambda^{-s}(-n_i{\rm div}u_i-u_i\cdot\nabla
 n_i)\Lambda^{-s}n_i
 -\Lambda^{-s}(u_i\cdot\nabla u_i+h(n_i)\nabla
 n_i)\cdot\Lambda^{-s}u_i\\[3mm]
 \leq &\D C\|n_i{\rm div}u_i+u_i\cdot\nabla
 n_i\|_{\dot{H}^{-s}}\|n_i\|_{\dot{H}^{-s}}+\|u_i\cdot\nabla
 u_i+h(n_i)\nabla n_i\|_{\dot{H}^{-s}}\|u_i\|_{\dot{H}^{-s}},
 \end{array}
 \eqno(2.56)
$$

If $s\in(0,1/2]$, then by Lemma 2.1, Lemma 2.3 and Young's
inequality, the right hand side of (2.56) can be estimated as
follows.
$$\arraycolsep=1.5pt  \begin{array}[b]{rl}
\|n_i{\rm div}u_i\|_{\dot{H}^{-s}} \leq &\D C\|n_i{\rm
div}u_i\|_{L^{\frac{1}{1/2+s/3}}}\leq C\|n_i\|_{L^{3/s}}\|\nabla
u_i\|_{L^2}\\[3mm]
\leq &\D C\|\nabla
n_i\|_{L^2}^{1/2+s}\|\nabla^2n_i\|_{L^2}^{1/2-s}\|\nabla
u_i\|_{L^2}\\[3mm]
\leq &\D C(\|\nabla n_i\|_{H^1}^2+\|\nabla u_i\|_{L^2}^2),
\end{array}
\eqno(2.57)
$$
where we have used the facts $\frac{1}{2}+\frac{s}{3}<1$ and
$\frac{3}{s}\geq6$.

Similarly, it holds that
$$
\|u_i\cdot\nabla n_i\|_{\dot{H}^{-s}}\leq C(\|\nabla
u_i\|_{H^1}^2+\|\nabla n_i\|_{L^2}^2);\eqno(2.58)
$$
$$
\|u_i\cdot\nabla u_i\|_{\dot{H}^{-s}}\leq C(\|\nabla
u_i\|_{H^1}^2+\|\nabla u_i\|_{L^2}^2);\eqno(2.59)
$$
$$
\|h(n_i)\cdot\nabla n_i\|_{\dot{H}^{-s}}\leq C(\|\nabla
n_i\|_{H^1}^2+\|\nabla n_i\|_{L^2}^2).\eqno(2.60)
$$

Now if $s\in(1/2,3/2)$, then $1/2+s/3<1$ and $2<3/s<6$. We shall
estimate the right hand side of (2.55) in a different way.  Using
Sobolev's inequality, we have
$$\arraycolsep=1.5pt  \begin{array}[b]{rl}
\|n_i{\rm div}u_i\|_{\dot{H}^{-s}} \leq &\D C\|n_i{\rm
div}u_i\|_{L^{\frac{1}{1/2+s/3}}}\leq C\|n_i\|_{L^{3/s}}\|\nabla
u_i\|_{L^2}\\[3mm]
\leq &\D C\|n_i\|_{L^2}^{s-1/2}\|\nabla n_i\|_{L^2}^{3/2-s}\|\nabla
u_i\|_{L^2},
\end{array}
\eqno(2.61)
$$
where we have used the facts $\frac{1}{2}+\frac{s}{3}<1$ and
$\frac{3}{s}\geq6$.

Similarly, it holds for $s\in(1/2,3/2)$ that
$$
\|u_i\cdot\nabla n_i\|_{\dot{H}^{-s}}\leq
C\|u_i\|_{L^2}^{s-1/2}\|\nabla u_i\|_{L^2}^{3/2-s}\|\nabla
n_i\|_{L^2};\eqno(2.62)
$$
$$
\|u_i\cdot\nabla u_i\|_{\dot{H}^{-s}}\leq
C\|u_i\|_{L^2}^{s-1/2}\|\nabla u_i\|_{L^2}^{3/2-s}\|\nabla
u_i\|_{L^2};\eqno(2.63)
$$
$$
\|h(n_i)\cdot\nabla n_i\|_{\dot{H}^{-s}}\leq
C\|n_i\|_{L^2}^{s-1/2}\|\nabla n_i\|_{L^2}^{3/2-s}\|\nabla
n_i\|_{L^2}.\eqno(2.64)
$$

Finally, we turn to the last term in the left hand side of (2.56)
with $i=1,2$. We have
$$\arraycolsep=1.5pt  \begin{array}[b]{rl}
&\D-\int\Lambda^{-s}\nabla\phi\cdot\Lambda^{-s}u_1+\int\Lambda^{-s}\nabla\phi\cdot\Lambda^{-s}u_2\\[3mm]
=&\D\int\Lambda^{-s}\phi\Lambda^{-s}{\rm
div}u_1-\int\Lambda^{-s}\phi\Lambda^{-s}{\rm
div}u_2\\[3mm]
=&\D-\int\Lambda^{-s}\phi\Lambda^{-s}\partial_t(n_1-n_2)+\int\Lambda^{-s}\phi\Lambda^{-s}{\rm
div}(n_1u_1-n_2u_2)\\[3mm]
=&\D\frac{1}{2}\frac{d}{dt}\int|\Lambda^{-s}\nabla\phi|^2-\int\Lambda^{-s}\nabla\phi\cdot\Lambda^{-s}(n_1u_1-n_2u_2).
\end{array}
\eqno(2.65)
$$

If $s\in(0,1/2)$, we use Lemma 2.1 and Lemma 2.4 to obtain
$$\arraycolsep=1.5pt  \begin{array}[b]{rl}
\|\Lambda^{-s}(n_iu_i)\|_{L^2}\leq &\D
C\|u_i\|_{L^2}\|n_i\|_{L^{3/s}}\leq C\|u_i\|_{L^2}\|\nabla
n_i\|_{L^2}^{1/2-s}\|\nabla^2n_i\|_{L^2}^{1/2+s}\\[3mm]
\leq &\D C(\|u_i\|_{L^2}^2+\|\nabla n_i\|_{H^1}^2);
\end{array}
\eqno(2.66)
$$
and if $s\in(1/2,3/2)$, we have
$$
\|\Lambda^{-s}(n_iu_i)\|_{L^2}\leq
C\|u_i\|_{L^2}\|n_i\|_{L^{3/s}}\leq C\|u_i\|_{L^2}\|\nabla
n_i\|_{L^2}^{s-1/2}\|\nabla^2n_i\|_{L^2}^{3/2-s}. \eqno(2.67)
$$

Consequently, in light of (2.56)-(2.67), and using Young's
inequality, we deduce (2.54) and (2.55). \qed

\subsection*{ 2.4.\ Estimates in $\dot{B}_{2,\infty}^{-s}(\mathbb{R}^3)$}

\quad\quad In this subsection, we will derive the evolution of the
negative Besov norms of the solutions. The argument is similar to
the previous subsection.

 \noindent\textbf{Lemma 2.14.} If
$\|n_{i0},u_{i0},\nabla\phi_0\|_{H^3}\ll1$ with $i=1,2$, for
$s\in(0,\frac{1}{2}]$, we have
$$
\frac{\rm d}{{\rm
d}t}\|(n_i,u_i,\nabla\phi)\|_{\dot{B}_{2,\infty}^{-s}}^2
 \leq
 C(\|\nabla n_i\|_{H^1}^2+\|u_i\|_{H^2}^2)\|(n_i,u_i,\nabla\phi)\|_{\dot{B}_{2,\infty}^{-s}},\
 i=1,2; \eqno(2.68)
$$
and for $s\in(\frac{1}{2},\frac{3}{2}]$, we have
$$\arraycolsep=1.5pt  \begin{array}[b]{rl}
 &\D \frac{\rm d}{{\rm d}t}\|(n_i,u_i,\nabla\phi)\|_{\dot{B}_{2,\infty}^{-s}}^2\\[3mm]
 \leq &\D
 C\left\{\|(n_i,u_i)|_{L^2}^{s-\frac{1}{2}}\|\nabla(n_i,u_i)\|_{H^1}^{\frac{5}{2}-s}
 +\|u_i\|_{L^2}\|n_i\|_{L^2}^{s-\frac{1}{2}}\|\nabla n_i\|_{L^2}^{\frac{3}{2}-s}\right\}\|(n_i,u_i,\nabla\phi)\|_{\dot{B}_{2,\infty}^{-s}},\
 i=1,2.
 \end{array}
 \eqno(2.69)
$$
\noindent\textbf{Proof.} Applying $\dot{\Delta}_j$ to $(2.2)_1$,
$(2.2)_2$ and multiplying the resulting identity by
$\dot{\Delta}_jn_1$, $\dot{\Delta}_ju_1$, respectively, and
integrating over $\mathbb{R}^3$ by parts, we get
$$\arraycolsep=1.5pt  \begin{array}[b]{rl}
 &\D \frac{\rm d}{{\rm d}t}\int(|\dot{\Delta}_jn_1|^2+|\dot{\Delta}_ju_1|^2)
 + \int|\nabla\dot{\Delta}_ju|^2-\int\dot{\Delta}_j\nabla\phi\cdot\dot{\Delta}_ju_1
\\[3mm]
 = &\D \int\dot{\Delta}_j(-n_1{\rm div}u_1-u_1\cdot\nabla n_1)\dot{\Delta}_jn_1
 -\dot{\Delta}_j(u_1\cdot\nabla u_1+h(n_1)\nabla
 n_1)\cdot\dot{\Delta}_ju_1\\[3mm]
 \leq &\D C\|n_1{\rm div}u_1+u_1\cdot\nabla
 n_1\|_{\dot{B}_{2,\infty}^{-s}}\|n_1\|_{\dot{B}_{2,\infty}^{-s}}+\|u_1\cdot\nabla
 u_1+h(n_1)\nabla
 n_1\|_{\dot{B}_{2,\infty}^{-s}}\|u_1\|_{\dot{B}_{2,\infty}^{-s}}.
 \end{array}
 \eqno(2.70)
$$

Then, as the proof of Lemma 2.13, applying Lemma 2.6 instead to
estimate the $\dot{B}_{2,\infty}^{-s}$ norm, we complete the proof
of Lemma 2.14. \qed

\section*{ 3.\ Proof of Theorems}

\subsection*{ 3.1.\ Proof of Theorem 1.1}

\quad\quad In this subsection, we shall use the energy estimates in
Subsection 2.2 to prove the global existence in $H^3$ norm.

We first close the energy estimates at each $l$-th level to prove
(1.3). Let $0\leq l\leq 2$. Summing up the estimates (2.4) from
$k=l$ to $k=2$, and then adding the resulting estimates to (2.31)
for $k=2$, by changing the index and since $\delta\ll1$, we have
$$\arraycolsep=1.5pt  \begin{array}[b]{rl}
&\D\frac{d}{dt}\sum_{l\leq k\leq
3}\|\nabla^k(n_1,u_1,n_2,u_2,\nabla\phi)\|_{L^2}^2+C_1\sum_{l\leq
k\leq 3}\|\nabla^k(u_1,u_2)\|_{L^2}^2\\[3mm]
\leq &\D C_2\delta\sum_{l+1\leq k\leq
3}\|\nabla^k(n_1,n_2,\nabla\phi)\|_{L^2}^2.
\end{array}
\eqno(3.1)
$$

Summing up (2.40) of Lemma 2.12 from $k=l$ to $2$, we have
$$
\arraycolsep=1.5pt  \begin{array}[b]{rl}
&\D\frac{d}{dt}\sum_{l\leq k\leq
2}\int(\nabla^ku_1\cdot\nabla\nabla^kn_1+\nabla^ku_2\cdot\nabla\nabla^kn_2)+C_3\sum_{l+1\leq
k\leq 3}\|\nabla^k(n_1,n_2,\nabla\phi)\|_{L^2}^2\\[3mm]
\leq &\D C_4\sum_{l\leq k\leq 3}\|\nabla^k(u_1,u_2)\|_{L^2}^2.
\end{array}
\eqno(3.2)
$$

Making a calculus $2C_2\delta/C_3\times(3.2)+(3.1)$, and by using
the fact $\delta\ll1$, we can conclude that there exists a constant
$C_5>0$ such that for $0\leq l\leq 2$,
$$
\arraycolsep=1.5pt  \begin{array}[b]{rl}
&\D\frac{d}{dt}\left\{\sum_{l\leq k\leq
3}\|\nabla^k(n_1,u_1,n_2,u_2,\nabla\phi)\|_{L^2}^2+\frac{2C_2\delta}{C_3}\sum_{l\leq
k\leq
m2}\int(\nabla^ku_1\cdot\nabla\nabla^kn_1+\nabla^ku_2\cdot\nabla\nabla^kn_2)\right\}\\[3mm]
&\D+C_5\left\{\sum_{l\leq k\leq
3}\|\nabla^k(u_1,u_2)\|_{L^2}^2+\sum_{l+1\leq k\leq
3}\|\nabla^k(n_1,n_2,\nabla\phi)\|_{L^2}^2\right\}\leq0.
\end{array}
\eqno(3.3)
$$
by the smallness of $\delta$ and using Cauchy's inequality, we
deduce that
$$
\arraycolsep=1.5pt  \begin{array}[b]{rl}
&\D C_6^{-1}\|\nabla^l(n_1,u_1,n_2,u_2,\nabla\phi)\|_{H^{3-l}}^2\\[3mm]
\leq &\D\sum_{l\leq k\leq
3}\|\nabla^k(n_1,u_1,n_2,u_2,\nabla\phi)\|_{L^2}^2+\frac{2C_2\delta}{C_3}\sum_{l\leq
k\leq
2}\int(\nabla^ku_1\cdot\nabla\nabla^kn_1+\nabla^ku_2\cdot\nabla\nabla^kn_2)\\[3mm]
\leq &\D C_6\|\nabla^l(n_1,u_1,n_2,u_2,\nabla\phi)\|_{H^{3-l}}^2,\
0\leq l\leq 2.
\end{array}
\eqno(3.4)
$$

As a result, we have the following estimate in Sobolev's space for
$0\leq l\leq 2$
$$
\frac{d}{dt}\|\nabla^l(n_1,u_1,n_2,u_2,\nabla\phi)\|_{H^{3-l}}^2
+\left\{\|\nabla^l(u_1,u_2)\|_{H^{3-l}}^2+\|\nabla^{l+1}(n_1,n_2,\nabla\phi)\|_{H^{2-l}}^2\right\}\leq0.
\eqno(3.5)
$$

Taking $l=0$ in (3.5), and integrating directly in time, we have
$$
\|(n_1,u_1,n_2,u_2,\nabla\phi)\|_{H^3}^2\leq
C_6^2\|(n_{10},u_{10},n_{20},u_{20},\nabla\phi_0)\|_{H^3}^2.\eqno(3.6)
$$

By a standard continuity argument, since
$\|(n_{10},u_{10},n_{20},u_{20},\nabla\phi_0)\|_{H^3}$ is
sufficiently small,  this closes the a priori estimates (2.2). Thus
we obtain the global existence in Theorem 1.1.

\subsection*{ 3.2.\ Proof of Theorem 1.2}

\quad\quad In this subsection, we will prove the optimal time decay
rates of the unique global solution to system (2.1) in theorem 1.2.

First, from Lemma 2.13, we need distinct the arguments by the value
of $s$. For $s\in[0,1/2)$, integrating (2.53) in time, by (1.3), we
have
$$\arraycolsep=1.5pt  \begin{array}[b]{rl}
\|(n_i,u_i,\nabla\phi)\|_{\dot{H}^{-s}}^2\leq &\D
\|(n_{i0},u_{i0},\nabla\phi_0)\|_{\dot{H}^{-s}}^2+C\int_0^t\|\nabla(n_i,u_i)\|_{H^1}^2\|(n_i,u_i,\nabla\phi)\|_{\dot{H}^{-s}}d\tau\\[3mm]
\leq &\D C_0(1+\sup\limits_{0\leq\tau\leq
t}\{\|(n_i,u_i,\nabla\phi)\|_{\dot{H}^{-s}}\}).
\end{array}
\eqno(3.7)
$$

This yields
$$
\|(n_1,u_1,n_2,u_2,\nabla\phi)\|_{\dot{H}^{-s}}\leq C_0\ {\rm for}\
s\in[0,1/2].\eqno(3.8)
$$

Using Lemma 2.14, we similarly have
$$
\|(n_1,u_1,n_2,u_2,\nabla\phi)\|_{\dot{B}_{2,\infty}^{-s}}\leq C_0\
{\rm for}\ s\in(0,1/2].\eqno(3.9)
$$

If $0\leq l\leq 2$, we may use Lemma 2.3 to have
$$
\|\nabla^{l+1}f\|_{L^2}\geq
C\|f\|_{\dot{H}^s}^{-\frac{1}{l+s}}\|\nabla^lf\|_{L^2}^{1+\frac{1}{l+s}}.\eqno(3.10)
$$

By this fact and (3.9), we find
$$
\|\nabla^{l+1}(n_1,n_2,\nabla\phi)\|_{L^2}^2\geq
C_0(\|\nabla^l(n_1,n_2,\nabla\phi)\|_{L^2}^2)^{1+\frac{1}{l+s}}.\eqno(3.11)
$$

This together with (1.3) yields for $l=0,1,2$,
$$
\|\nabla^l(u_1,u_2),\nabla^{l+1}(n_1,n_2,\nabla\phi)\|_{H^{3-l}}^2\geq
C_0(\|\nabla^l(u_1,u_2,n_1,n_2,\nabla\phi)\|_{H^{3-l}}^2)^{1+\frac{1}{l+s}}.\eqno(3.12)
$$

Hence, from (3.5), we have the following time differential
inequality for $l=0,1,2$
$$
\frac{d}{dt}\|\nabla^l(u_1,u_2,n_1,n_2,\nabla\phi)\|_{H^{3-l}}^2
+C_0(\|\nabla^l(u_1,u_2,n_1,n_2,\nabla\phi)\|_{H^{3-l}}^2)^{1+\frac{1}{l+s}}\leq0,\eqno(3.13)
$$
which gives
$$
\|\nabla^l(u_1,u_2,n_1,n_2,\nabla\phi)\|_{H^{3-l}}^2\leq
C_0(1+t)^{-(l+s)},\ l=0,1,2;\ s\in[0,\frac{1}{2}].\eqno(3.14)
$$

For $s\in(1/2,3/2)$. Notice that the arguments for the case
$s\in[0,1/2]$ can not be applied to this case (See Lemma 2.13).
Observing that we have
$n_{10},u_{10},n_{20},u_{20},\nabla\phi_0\in\dot{H}^{-1/2}$ since
$\dot{H}^{-s}\cap L^2\subset \dot{H}^{-s'}$ for any $s'\in[0,s]$, we
then deduce from what we have proved for (1.6) with $s=1/2$ that the
following decay result holds:
$$
\|\nabla^l(n_1,u_1,n_2,u_2,\nabla\phi)\|_{H^{3-l}}\leq
C_0(1+t)^{-\frac{l+\frac{1}{2}}{2}}\ {\rm for}\ l=0,1,2.\eqno(3.15)
$$

Integrating (2.53) in time, for $s\in(1/2,3/2)$, we have
$$\arraycolsep=1.5pt  \begin{array}[b]{rl}
&\D\|(n_i,u_i,\nabla\phi)\|_{\dot{H}^{-s}} \leq
\|(n_{i0},u_{i0},u_{20},\nabla\phi_0)\|_{\dot{H}^{-s}}\\[3mm]
&\D+C\int_0^t\left\{\|(n_i,u_i)\|_{L^2}^{s-\frac{1}{2}}\|\nabla(n_i,u_i)\|_{H^1}^{\frac{5}{2}-s}
 +\|u_i\|_{L^2}\|n_i\|_{L^2}^{s-\frac{1}{2}}\|\nabla
 n_i\|_{L^2}^{\frac{3}{2}-s}\right\}\|(n_i,u_i,\nabla\phi)\|_{\dot{H}^{-s}}d\tau\\[3mm]
 \leq &\D\|(n_{i0},u_{i0},u_{20},\nabla\phi_0)\|_{\dot{H}^{-s}}
+C\sup_{0\leq\tau\leq
t}\{\|(n_i,u_i,\nabla\phi)\|_{\dot{H}^{-s}}\}\\[3mm]
&\D\ \ \ \ \ \ \
 \ \ \  \ \ \ \ \ \ \ \ \ \ \ \ \ \ \times\int_0^t\left\{\|(n_i,u_i)\|_{L^2}^{s-\frac{1}{2}}\|\nabla(n_i,u_i)\|_{H^1}^{\frac{5}{2}-s}
 +\|u_i\|_{L^2}\|n_i\|_{L^2}^{s-\frac{1}{2}}\|\nabla
 n_i\|_{L^2}^{\frac{3}{2}-s}\right\}d\tau\\[3mm]
 :=&\D\|(n_{i0},u_{i0},u_{20},\nabla\phi_0)\|_{\dot{H}^{-s}}+C\sup_{0\leq\tau\leq
t}\{\|(n_i,u_i,\nabla\phi)\|_{\dot{H}^{-s}}\}\cdot(K_1+K_2).
 \end{array}
 \eqno(3.16)
$$

For $K_1$, by using (3.15), we deduce that for the case
$s\in(\frac{1}{2},\frac{3}{2})$
$$
\arraycolsep=1.5pt  \begin{array}[b]{rl} K_1=&\D
C\int_0^t\{\|(n_i,u_i)\|_{L^2}^{s-\frac{1}{2}}\|\nabla(n_i,u_i)\|_{H^1}^{\frac{5}{2}-s}
 \}\|(n_i,u_i,\nabla\phi)\|_{\dot{H}^{-s}}d\tau\\[3mm]
\leq &\D C_0+C_0\int_0^t(1+\tau)^{-7/4-s/2}d\tau\sup_{0\leq\tau\leq
t}\{\|(n_i,u_i,\nabla\phi)\|_{\dot{H}^{-s}}\}\\[3mm]
\leq &\D C_0\left\{1+\sup_{0\leq\tau\leq
t}\{\|(n_i,u_i,\nabla\phi)\|_{\dot{H}^{-s}}\})\right\},\ i=1,2;\
s\in(\frac{1}{2},\frac{3}{2}).
\end{array}
\eqno(3.17)
$$

For $K_2$, we must distinct the arguments by the value of $s$:
$s\in(\frac{1}{2},1)$ and $s\in[1,\frac{3}{2})$. When
$s\in(\frac{1}{2},1)$,
$$\arraycolsep=1.5pt  \begin{array}[b]{rl}
K_2=&\D
C\int_0^t\{\|u_i\|_{L^2}\|n_i\|_{L^2}^{s-\frac{1}{2}}\|\nabla
 n_i\|_{L^2}^{\frac{3}{2}-s}\}\|(n_i,u_i,\nabla\phi)\|_{\dot{H}^{-s}}d\tau\\[3mm]
 \leq &\D C\{\int_0^t\|u_i\|_{L^2}^2d\tau+\int_0^t\|n_i\|_{L^2}^{2s-1}\|\nabla
 n_i\|_{L^2}^{3-2s}d\tau\\[3mm]
 \leq &\D
 CC_0+C_0\int_0^t(1+\tau)^{-\frac{1}{4}(2s-1)}(1+\tau)^{-\frac{3}{4}(3-2s)}d\tau\\[3mm]
 \leq &\D CC_0+\int_0^t(1+\tau)^{-\frac{1}{4}(8-4s)}d\tau\leq CC_0,\
 s\in(\frac{1}{2},1).
 \end{array}
 \eqno(3.18)
 $$

Thus, (3.16)-(3.18) imply that
$$
\|(n_i,u_i,\nabla\phi)\|_{\dot{H}^{-s}}\leq CC_0,\ s\in[0,1).
\eqno(3.19)
$$

Hence (3.19) together with a similar argument as the case
$s\in[0,\frac{1}{2}]$, we know the decay result (1.6) is established
for any $s\in[0,1)$:
$$
\|\nabla^l(u_1,u_2,n_1,n_2,\nabla\phi)\|_{H^{3-l}}^2\leq
C_0(1+t)^{-(l+s)},\ l=0,1,2;\ s\in[0,1).\eqno(3.20)
$$

Now we choose a constant $s_1=\frac{5}{8}+\frac{s}{4}$ with
$s\in[1,\frac{3}{2})$, then $s_1<1$. Then, (3.20) gives
$$
\|\nabla^l(u_1,u_2,n_1,n_2,\nabla\phi)\|_{H^{3-l}}^2\leq
C_0(1+t)^{-(l+s_1)},\ l=0,1,2;\ s_1\in[0,1).\eqno(3.21)
$$

By (3.21), we can prove the decay result for $s\in[1,\frac{3}{2})$.
In fact,
$$\arraycolsep=1.5pt  \begin{array}[b]{rl}
K_2=&\D
C\int_0^t\{\|u_i\|_{L^2}\|n_i\|_{L^2}^{s-\frac{1}{2}}\|\nabla
 n_i\|_{L^2}^{\frac{3}{2}-s}\}\|(n_i,u_i,\nabla\phi)\|_{\dot{H}^{-s}}d\tau\\[3mm]
\leq &\D
CC_0\int_0^t(1+\tau)^{-\frac{s_1}{2}}(1+\tau)^{-\frac{s_1}{2}(s-\frac{1}{2})}(1+\tau)^{-\frac{1+s_1}{2}(\frac{3}{2}-s)}d\tau\\[3mm]
=&\D CC_0\int_0^t(1+\tau)^{s_1+\frac{3}{4}-\frac{s}{2}}d\tau=
CC_0\int_0^t(1+\tau)^{\frac{11}{8}-\frac{s}{4}}d\tau\leq CC_0,\
s\in[1,\frac{3}{2}).
\end{array}
\eqno(3.22)
$$

Hence, (3.16), (3.17) and (3.22) suffice for that
$$
\|(n_i,u_i,\nabla\phi)\|_{\dot{H}^{-s}}\leq CC_0,\
s\in[0,\frac{3}{2}). \eqno(3.23)
$$

With (3.23) in hand, we repeat the arguments leading to (1.6) for
$s\in[0,1/2]$ to prove that it hold also for $s\in(1/2.3/2)$.

Lastly, by using Lemma 2.5, Lemma 2.7, Lemma 2.8, Lemma 2.9 and
Lemma 2.14, a similar argument as leading to the estimate (3.23) for
the negative Sobolev space can immediately yields that in the
negative Besov's space:
$$
\|(n_i,u_i,\nabla\phi)\|_{\dot{B}_{2,\infty}^{-s}}\leq CC_0,\
s\in(0,\frac{3}{2}]. \eqno(3.24)
$$
\\
\\
{\bf Acknowledgement:} \ \  The research of Zhigang Wu  was
supported by NSFC (No. 11101112) and in part by NSFC (No. 11271105)
and the Nature Science Foundation of Zhejiang Province (No.
LY12A01030). The research of W.K. Wang was supported by the National
Science Foundation of China 11071162.

\bibliographystyle{plain}

\begin{thebibliography}{99}

\bibitem{Ali}G. Ali, A. J\"{u}gel, Global smooth solutions to the
multi-dimensional hydrodynamic model for two-carrier plasma,
\emph{J. Differential Equations} \textbf{190}(2003), 663-685.

\bibitem{Gasser1}I. Gasser, L. Hsiao, H.L. Li, Asymptotic behavior of solutions
of the bipolar hydrodynamic fluids, \emph{J. Differential Equations}
\textbf{192}(2003), 326-359.

\bibitem{Gasser2}I. Gasser, P. Marcati, The combined relaxation and vanishing
Debye length limit in the hydrodynamic model for semiconductors,
\emph{Math. Methods Appl. Sci.} \textbf{24}(2001), 81-92.

\bibitem{Guo}Y. Guo, Y.J. Wang, Decay of dissipative equations and negative Sobolev spaces,
\emph{Comm. Partial Differential Equations} \textbf{37}(2012),
2165-2208.

\bibitem{Jugel}A. J\"{u}ngel, Quasi-Hydrodynamic Semiconductor Equations, \emph{Progr.
Nonlinear Differential Equations Appl.}, Birkhauser, 2001.

\bibitem{Ju}Q.C. Ju, Global smooth solutions to the multidimensional
hydrodynamic model for plasmas with insulating boundary conditions,
\emph{J. Math. Anal. Appl.} \textbf{336}(2007), 888-904.

\bibitem{Huang}F.M. Huang, Y.P. Li, Large time behavior and quasineutral limit
of solutions to a bipolar hydrodynamic model with large data and
vacuum, \emph{Discrete Contin. Dyn. Syst. Ser. A} \textbf{24}(2009),
455-470.

\bibitem{Hsiao1}L. Hsiao, K.J. Zhang, The global weak solution and relaxation
limits of the initial boundary value problem to the bipolar
hydrodynamic model for semiconductors, \emph{Math. Models Methods
Appl. Sci.} \textbf{10}(2000), 1333-1361.

\bibitem{Hsiao2}L. Hsiao, K.J. Zhang, The relaxation of the hydrodynamic model
for semiconductors to the drift-diffusion equations, \emph{J.
Differential Equations} \textbf{165}(2000) 315-354.

\bibitem{Lattanzio}C. Lattanzio, On the 3-D bipolar isentropic Euler-Poisson model
for semiconductors and the drift-diffusion limit, \emph{Math. Models
Methods Appl. Sci.} \textbf{10}(2000), 351-360.

\bibitem{LiH1}H.L. Li, A. Matsumura, G.J. Zhang, Optimal decay rate of the
compressible Navier-Stokes-Poisson system in $R^3$, Arch Ration Mech
Anal, 2010, \emph{196}: 681-713.

\bibitem{LiH2}H.L. Li, T. Yang, C. Zou, Time asymptotic behavior of the
bipolar Navier-Stokes-Poisson system, \emph{Acta. Mathematica
Scientia}, 2009, \textbf{29B}(6): 1721-1736.

\bibitem{Li1}Y.P. Li, T. Zhang, Relaxation-time limit of the multidimensional
bipolar hydrodynamic model in Besov space, \emph{J. Differential
Equations} \textbf{251}(11)(2011), 3143-3162.

\bibitem{Li2}Y.P. Li, Diffusion relaxation limit of a bipolar isentropic
hydrodynamic model for semiconductors, \emph{J. Math. Anal. Appl.}
\textbf{336}(2007), 1341-1356.

\bibitem{Li3}Y.P. Li, X.F. Yang, Global existence and asymptotic behavior of the solutions
to the three-dimensional bipolar Euler-Poisson systems, \emph{J.
Differential Equations} \textbf{252}(2012), 768-791.

\bibitem{Kato}T. Kato, The Cauchy problem for quasi-linear symmetric
hyperbolic systems, \emph{Arch. Ration. Mech. Anal.} \textbf{58}
(1945), 181-205.

\bibitem{Markowich}P.A. Markowich, C.A. Ringhofev, C. Schmeiser, Semiconductor
Equations, Springer, Wien, New York, 1990.

\bibitem{Natalini1}R. Natalini, The bipolar hydrodynamic model for semiconductors
and the drift-diffusion equation, \emph{J. Math. Anal. Appl.}
\textbf{198}(1996), 262-281.

\bibitem{Sitenko}A. Sitenko, V. Malnev, Plasma Physics Theory, Chapman \& Hall,
London, 1995.

\bibitem{Sohinger}V. Sohinger, R.M. Strain, The Boltzmann equation, Besov spaces,
and optimal time decay rates in Rnx , preprint, arXiv: 1206.0027,
2012.

\bibitem{Stein}E. M. Stein, Singular Integrals and Differentiability Properties
of Functions, Princeton University Press, 1970.

\bibitem{Tan}Z. Tan, Y. Wang, Global solution and large-time behavior of the
3D compressible Euler equations with damping, \emph{J. Differential
Equations}, in press.

\bibitem{Taylor}M.E. Taylor, Partial Differential Equations (III): Nonlinear
Equations, Springer, 1996.

\bibitem{Tsuge}N. Tsuge, Existence and uniqueness of stationary solutions to a
one-dimensional bipolar hydrodynamic models of semiconductors,
\emph{Nonlinear Anal.} \textbf{73}(2010), 779-787.

\bibitem{Wang1}W.K. Wang, Z.G. Wu, Pointwise estimates of solution for
the Navier-Stokes-Poisson equations in multi-dimensions, \emph{J.
Differential Equations} \textbf{248}(2010), 1617-1636.

\bibitem{WangY1}Y.J. Wang, Decay of the Navier-Stokes-Poisson
equations, \emph{J. Differential Equations} \textbf{253}(2012),
273-297.

\bibitem{Wu1}Z.G. Wu, W.K. Wang, Pointwise estimates of solution for
non-isentropic Navier-Stokes equations in multi-dimensions,
\emph{Acta Mathematica Scientia} \textbf{32B}(5)(2012), 1681-1702.

\bibitem{Wu2}Z.G. Wu, W.K. Wang, Pointwise estimate of solutions for the Euler-Poisson equation
 with damping in multi-dimensions, \emph{Discrete Contin. Dyn. Syst. A},
\textbf{26}(3)(2010), 1101-1117.

\bibitem{Wu3}Z.G. Wu, A refined energy estimate for the Euler-Poisson equation
 with damping in 3D, submitted.

\bibitem{Zhang}G.J. Zhang, H.L. Li, C.J. Zhu, Optimal decay rate of the
non-isentropic Navier-Stokes-Poisson system in $R^3$, \emph{J.
Differential Equations} \textbf{250}(2)(2011), 866-891.

\bibitem{Zhou}F. Zhou, Y.P. Li, Existence and some limits of stationary
solutions to a one-dimensional bipolar Euler¨CPoisson system,
\emph{J. Math. Anal. Appl.} \textbf{351}(2009), 480-490.

\bibitem{Zhu}C. Zhu, H. Hattori, Stability of steady state solutions for an
isentropic hydrodynamic model of semiconductors of two species,
\emph{J. Differential Equations} \textbf{166}(2000), 1-32.






\end{thebibliography}

\end{document}